\documentclass[sn-mathphys-num]{sn-jnl}

\usepackage{graphicx}%
\usepackage{multirow}%
\usepackage{amsmath,amssymb,amsfonts}%
\usepackage{amsthm}%
\usepackage{mathrsfs}%
\usepackage[title]{appendix}%
\usepackage{xcolor}%
\usepackage{textcomp}%
\usepackage{manyfoot}%
\usepackage{booktabs}%
\usepackage{algorithm}%
\usepackage{algorithmicx}%
\usepackage{algpseudocode}%
\usepackage{listings}%

\theoremstyle{thmstyleone}%
\newtheorem{theorem}{Theorem}[section]
\newtheorem{proposition}[theorem]{Proposition}%

\theoremstyle{thmstyletwo}%
\newtheorem{remark}[theorem]{Remark}%
\newtheorem{lemma}[theorem]{Lemma}
\newtheorem{corollary}[theorem]{Corollary}
\theoremstyle{thmstylethree}%
\newtheorem{definition}[theorem]{Definition}%
\begin{document}

\title[Article Title]{Counting of lattices containing up to four comparable reducible elements and having nullity up to three}


\author*[1,2]{\fnm{B. P.} \sur{Aware}}\email{aware66pp@gmail.com}

\author[2,3]{\fnm{Dr. A. N.} \sur{Bhavale}}\email{ hodmaths@moderncollegepune.edu.in}
\equalcont{These authors contributed equally to this work.}

\affil*[1]{\orgdiv{Department of Mathematics}, \orgname{Modern College of ASC(A)}, \orgaddress{\street{Shivajinagar}, \city{Pune}, \postcode{411005}, \state{Maharastra}, \country{India}}}

\affil[2]{\orgdiv{Head, Department of Mathematics}, \orgname{Modern College of ASC(A)}, \orgaddress{\street{Shivajinagar}, \city{Pune}, \postcode{411005}, \state{Maharastra}, \country{India}}}

\abstract{In 2020 Bhavale and Waphare introduced the concept of a nullity of a poset as nullity of its cover graph. According to Bhavale and Waphare, if a dismantlable lattice of nullity k contains r reducible elements then 2 $\leq$ r $\leq$ 2k. In 2003 Pawar and Waphare counted all non-isomorphic lattices with equal number of elements and edges, which are precisely the lattices of nullity one. Recently, Bhavale and Aware counted all non-isomorphic lattices on n elements having nullity up to two. Bhavale and Aware also counted all non-isomorphic lattices on n elements, containing up to three reducible elements, having nullity k $\geq$ 2. In this paper, we count up to isomorphism the class of all lattices on n elements containing four comparable reducible elements, and having nullity three.}



\maketitle

\section{Introduction}\label{sec1}

Birkhoff \cite{bib6} in $1940$ raised the open problem to compute for small $n$ all posets/lattices on a set of $n$ elements up to isomorphism. Many authors all over the world attempted this problem. Enumeration of all posets with $15$ and $16$ elements up to isomorphism was done by Brinkmann and Mckay \cite{bib7}. Enumeration of all posets up to isomorphism is still in progress for $n\geq 17$. Counting of all lattices on up to $18$ elements up to isomorphism was carried out by Heitzig and Reinhold \cite{bib9}.

Bhavale and Waphare \cite{bib5} defined the concept of a nullity of a poset $P$ denoted by $\eta(P)$ as nullity of its cover graph which is given by $p-q+c$, where $p$ is the number of edges, $q$ is the number of vertices, and $c$ is the number of connected components of the graph. In the same paper, they also defined the concept of an RC-lattice as a lattice in which all the reducible elements are comparable. Bhavale and Waphare \cite{bib5}, shown that for a dismantlable lattice having nullity $k$ containing $r$ reducible elements satisfies the inequality $2\leq r\leq 2k$. In $2003$ Pawar and Waphare \cite{bib11} counted all non-isomorphic lattices with $n$ elements, containing exactly two reducible elements, and having nullity one. Thakare et al. \cite{bib14} generalised this result for arbitrary nullity $k\geq 1$. Recently, Bhavale and Aware \cite{bib2} counted all non-isomorphic lattices on $n$ elements, containing up to three reducible elements, and having arbitrary nullity $k\geq 2$. Also, Bhavale and Aware \cite{bib3} counted all non-isomorphic lattices on $n$ elements having nullity up to two. Thus, this counting covers counting of all non-isomorphic lattices on $n$ elements containing $r\leq 4$ reducible elements, and having nullity $k\leq 2$(see \cite{bib2}, \cite{bib3}). That means, the problem of counting of all non-isomorphic lattices on $n$ elements containing $r\geq 4$ reducible elements, and having nullity $k\geq 3$ is still open. In this paper, we count up to isomorphism the class of all non-isomorphic RC-lattices on $n$ elements, containing exactly four reducible elements, and having nullity three.

An element $y$ in $P$ {\it{covers}} an element $x$ in $P$ if $x<y$, and there is no element $z$ in $P$ such that $x<z<y$. We denote this fact by $x\prec y$, and say that pair $<x,y>$ is a $covering$ or an $edge$({\it{see \cite{bib14}}}). An element of a poset $P$ having at most one lower cover and at most one upper cover in $P$ is called {\it{doubly irreducible}}. The set of all doubly irreducible elements in the poset $P$ is denoted by $Irr(P)$. The set of all doubly irreducible elements having exactly one upper cover and exactly one lower cover in $P$ is denoted by $Irr^*(P)$. The graph on a poset $P$ with edges as covering relations is called the {\it{cover graph}} and is denoted by $C(P)$. The number of coverings in a chain $C$ is called {\it{length}}, denoted by $|C|$ of that chain.

 An element $a$ in a lattice $L$ is $meet$-$reducible$($join$-$reducible$) in $L$ if there exist $b,c \in L$ both distinct from $a$, such that $inf\{b,c\}=a (sup\{b,c\}=a)$. An element $a$ in a lattice $L$ is said to be $reducible$ if it is either meet-reducible or join-reducible. $a$ is said to be $meet$-$irreducible$($join$-$irreducible$) if it is not meet-reducible(join-reducible). $a$ is $doubly$ $irreducible$ if it is both meet-irreducible and join-irreducible. The set of all doubly irreducible elements in $L$ is denoted by $Irr(L)$, and $Red(L)$ denotes  complement of $Irr(L)$ in $L$.

The following definition is due to 1974 Rival \cite{bib12}. 
\begin{definition}\cite{bib12}
\textnormal{A finite lattice of order $n$ is called {\it{dismantlable}} if there exists a chain $L_{1} \subset L_{2} \subset \ldots\subset L_{n}(=L)$ of sublattices of $L$ such that $|L_{i}| = i$, for all $i$}.
\end{definition}

The concept of {\it{adjunct operation of lattices}}, is introduced by Thakare et al. \cite{bib14}. 
Suppose $L_1$ and $L_2$ are two disjoint lattices and $(a, b)$ is a pair of elements in $L_1$ such that $a<b$ and $a\not\prec b$. Define the partial order $\leq$ on $L = L_1 \cup L_2$ with respect to the pair $(a,b)$ as follows: $x \leq y$ in L if $x,y \in L_1$ and $x \leq y$ in $L_1$, or $x,y \in L_2$ and $x \leq y$ in $L_2$, or $x \in L_1,$ $ y \in L_2$ and $x \leq a$ in $L_1$, or $x \in L_2,$ $ y \in L_1$ and $b \leq y$ in $L_1$. 

It is easy to see that $L$ is a lattice containing $L_1$ and $L_2$ as sublattices. The procedure for obtaining $L$ in this way is called {\it{adjunct operation (or adjunct sum)}} of $L_1$ with $L_2$. We call the pair $(a,b)$ as an {\it{adjunct pair}} and $L$ as an {\it{adjunct}} of $L_1$ with $L_2$ with respect to the adjunct pair $(a,b)$ and write $L = L_1 ]^b_a L_2$. A diagram of $L$ is obtained by placing a diagram of $L_1$ and a diagram of $L_2$ side by side in such a way that the largest element $1$ of $L_2$ is at lower position than $b$ and the least element $0$ of $L_2$ is at the higher position than $a$ and then by adding the coverings $<1,b>$ and $<a,0>$. This clearly gives $|E(L)|=|E(L_1)|+|E(L_2)|+2$. A lattice $L$ is called an {\it{adjunct of lattices}} $L_0,L_1,\ldots,L_k,$ if it is of the form $L = L_0]^{b_1}_{a_1} L_1 ]^{b_2}_{a_2} L_2 \cdots ]^{b_{k}}_{a_{k}} L_k$.
\begin{theorem}\cite{bib14}\label{dac} 
A finite lattice is dismantlable if and only if it is an adjunct of chains.
\end{theorem}
If $P_1$ and $P_2$ are two disjoint posets, the {\it{direct sum}} ({\it{see \cite{bib13}}}), denoted by $P_1 \oplus P_2$, is defined by taking the following order relation on $P_1\cup P_2:x\leq y$ if and only if $x,y\in P_1$ and $x\leq y$ in $P_1$, or $x,y\in P_2$ and $x\leq y$ in $P_2$ or $x\in P_1,y\in P_2$. In general $P_1 \oplus P_2\neq P_2 \oplus P_1$. Also, if $P_1$ and $P_2$ are lattices then $|E(P_1\oplus P_2)|=|E(P_1)|+|E(P_2)|+1$.

\begin{lemma}\label{oplus}
Let $\mathscr{L}_{1}(p)$ and $\mathscr{L}_{2}(q)$ be some classes of non-isomorphic lattices on $p\geq r\geq 1$ and $q\geq s\geq 1$ elements respectively. If $L_1\in \mathscr{L}_{1}(p)$ and $L_2\in \mathscr{L}_{2}(q)$ are such that $L_1\oplus L_2=L\in \mathscr{L}(n)$, a class of non-isomorphic lattices on $n=p+q$ elements then $|\mathscr{L}(n)|=\displaystyle\sum_{p=r}^{n-s}(|\mathscr{L}_{1}(p)|\times|\mathscr{L}_{2}(n-p)|)$.
\end{lemma}
\begin{proof}
As $n=p+q$, for fixed $p$, there are $|\mathscr{L}_{1}(p)|$ non-isomorphic lattices on $p$ elements and $|\mathscr{L}_{2}(n-p)|$ non-isomorphic lattices on $n-p$ elements. Therefore by multiplication principle, for fixed $p$, there are up to isomorphism $|\mathscr{L}_{1}(p)|\times|\mathscr{L}_{2}(n-p)|$ lattices in $\mathscr{L}(n)$. Note that $r\leq p=n-q\leq n-s$, since $q\geq s$. Thus there are $\displaystyle\sum_{p=r}^{n-s}(|\mathscr{L}_{1}(p)|\times|\mathscr{L}_{2}(n-p)|)$ lattices in $\mathscr{L}(n)$ up to isomorphism.
\end{proof}
 The \textit{vertical sum} of posets $P_1$ and $P_2$, is obtained from the the {\it{direct sum}} $P_1 \oplus P_2$ by identifying greatest element of $P_1$ with the least element of $P_2$ (see \cite{bib16}). The following definition is due to Bhamre and Pawar \cite{bib4}.
\begin{definition}\cite{bib4}
\textnormal{Let $P_1$ be a poset with the largest element and $P_2$ be a poset with the least element such that the greatest element of $P_1$ and the least element of $P_2$ are the same(say $\alpha$) and $P_1\cap P_2 =\{\alpha\}$, then the vertical sum of $P_1$ with $P_2$, denoted by $P_1 \circ P_2$, is a poset $(P_1\cup P_2, \leq)$ where $x \leq y$ if and only if $x, y \in P_1$ and $x \leq y$ in $P_1$, or $x, y\in P_2$ and $x \leq y$ in $P_2$, or $x$ in $P_1$ and $y$ in $P_2$.}
\end{definition}
\begin{lemma}\label{circ}
Let $\mathscr{L}_{1}(p)$ and $\mathscr{L}_{2}(q)$ be some classes of non-isomorphic lattices on $p\geq r\geq 1$ and $q\geq s\geq 1$ elements respectively. If $L_1\in \mathscr{L}_{1}(p)$ and $L_2\in \mathscr{L}_{2}(q)$ are such that $L_1\circ L_2=L\in \mathscr{L}(n)$, a class of non-isomorphic lattices on $n=p+q-1$ elements then $|\mathscr{L}(n)|=\displaystyle\sum_{p=r}^{n-s+1}(|\mathscr{L}_{1}(p)|\times|\mathscr{L}_{2}(n-p+1)|)$.
\end{lemma}
\begin{proof}
As $n=p+q-1$, for fixed $p$, there are $|\mathscr{L}_{1}(p)|$ non-isomorphic lattices. For $m\in\mathbb{N}$, let $\mathscr{P}_2(m)=\{L\setminus\{0\}~|~ L\in\mathscr{L}_{2}(m+1)\}$. Observe that there is a one-one correspondence between the class $\mathscr{P}_2(m)$ and the class $\mathscr{L}_2(m+1)$, and hence $|\mathscr{L}_2(m+1)|=|\mathscr{P}_2(m)|$. Note that for fixed $p$, there are $|\mathscr{P}(n-p)|$ non-isomorphic posets on $n-p$ elements, and hence there are $|\mathscr{L}_{2}(n-p+1)|$ non-isomorphic lattices on $n-p+1$ elements. Therefore by multiplication principle, for fixed $p$, there are $|\mathscr{L}_{1}(p)|\times|\mathscr{L}_{2}(n-p+1)|$ lattices on $n$ elements in $\mathscr{L}(n)$. Note that $1 \leq r\leq p=n-q+1\leq n-s+1$, since $q\geq s\geq 1$. Thus there are $\displaystyle\sum_{p=r}^{n-s+1}(|\mathscr{L}_{1}(p)|\times|\mathscr{L}_{2}(n-p+1)|)$ lattices in $\mathscr{L}(n)$ up to isomorphism.
\end{proof}
Thakare et al. \cite{bib14} defined a {\it block} as a finite lattice in which the largest element 1 is join-reducible and the smallest element 0 is meet-reducible. Moreover, if $L$ is any lattice which is different from chain, then $L$ contains a unique maximal sublattice which is a block called as {\it{maximal block}}. The lattice $L$ is of the form $C\oplus\textbf{B}$ or $\textbf{B}\oplus C$ or $C\oplus\textbf{B}\oplus C'$, where $C,C'$ are chains and $\textbf{B}$ is the maximal block. Further $\eta(L)=\eta(\textbf{B})$.
Bhavale and Waphare \cite{bib5} introduced the following concepts namely, retractible element, basic retract, basic block, basic block associated to a poset.
\begin{definition}\cite{bib5}\label{rtrct}
\textnormal{Let $P$ be a poset. Let $x \in Irr(P)$. Then $x$ is called a {\it{retractible}} element of $P$ if it satisfies either of the following conditions.}
\begin{enumerate}
\item \textnormal{There are no $y,z \in Red(P)$ such that $y \prec x \prec z$}.
\item \textnormal{There are $y,z \in Red(P)$ such that $y \prec x \prec z$ and there is no other directed path from $y$ to $z$ in $P$.}
\end{enumerate}
\end{definition}
\begin{definition}\cite{bib5} \label{brt}
\textnormal{A poset $P$ is a $basic$ $retract$ if no element of $Irr^{*}(P)$ is retractible in the poset $P$}.
\end{definition}
\begin{definition}\cite{bib5}\label{basicblock}
\textnormal{A poset $P$ is a {\it{basic block}} if it is one element or $Irr(P) = \phi$ or removal of any doubly irreducible element reduces nullity by one}.
\end{definition}
\begin{definition}\cite{bib5}\label{bbas}
\textnormal{$B$ is a $basic$ $block$ $associated$ $to$ $a$ $poset$ $P$ if $B$ is obtained from the basic retract associated to $P$ by successive removal of all the pendant vertices.}
\end{definition}
\begin{theorem}\cite{bib5}\label{redb}
Let $B$ be the basic block associated to a poset $P$. Then $Red(B)=Red(P)$ and $\eta(B)=\eta(P)$.
\end{theorem}
\begin{theorem}\cite{bib5}\label{nulladj}
A dismantlable lattice $L$ containing $n$ elements is of nullity $l$ if and only if $L$ is an adjunct of $l+1$ chains.
\end{theorem}
For the other definitions, notation, and terminology, see \cite{bib2}, \cite{bib8} and \cite{bib15}.

\section{Counting of lattices containing up to three reducible elements}\label{sec2}
In this section, we discuss counting of all non-isomorphic lattices on $n$ elements containing up to three reducible elements. 
 Let $\mathscr{L}(n;r,k)$ be the class of all non-isomorphic lattices on $n$ elements such that every member of it contains $r$ comparable reducible elements, and has nullity $k$. Let $\mathscr{L}(n;r,k,h)$ be the subclass of $\mathscr{L}(n;r,k)$ such that the basic block associated to a member of it is of height $h$. Let $\mathscr{B}(j;r,k)$ be the class of all non-isomorphic maximal blocks on $j$ elements such that every member of it contains $r$ comparable reducible elements and has nullity $k$. Let $\mathscr{B}(j;r,k,h)$ be the subclass of $\mathscr{B}(j;r,k)$ such that the basic block associated to a member of it is of height $h$.  Let $P_n^k$ denote the number of partitions of $n$ into $k$ non-decreasing positive integer parts.

 It is clear that chain is the only lattice containing no reducible element and having nullity $0$. Obviously, there is no lattice on one reducible element. Thakare et. al. \cite{bib14} enumerated all non-isomorphic lattices on $n$ elements, containing exactly two reducible elements, and having nullity $k$.
\begin{theorem}(\cite{bib14},\cite{bib2})\label{2red}
\textnormal{For} $n\geq 4$ and for $1\leq k\leq n-3$, $|\mathscr{L}(n;2,k)|=\displaystyle\sum_{j=1}^{n-k-2}jP_{n-j-1}^{k+1}$.
\end{theorem}

\begin{center}
\unitlength 1mm 
\linethickness{0.4pt}
\ifx\plotpoint\undefined\newsavebox{\plotpoint}\fi 
\begin{picture}(123.268,37.054)(0,0)
\put(9.362,30.213){\circle{2.154}}
\put(9.362,22.476){\circle{2.154}}
\put(9.53,15.581){\circle{2.154}}
\put(1.122,26.008){\circle{2.154}}
\put(9.362,35.594){\circle{2.154}}
\put(9.126,17.027){\line(0,1){0}}
\put(9.126,17.027){\line(0,1){0}}
\put(9.378,34.434){\line(0,-1){3.027}}
\put(9.378,31.406){\line(0,1){0}}
\put(9.294,29.136){\line(0,-1){5.634}}
\put(9.378,21.4){\line(0,-1){4.793}}
\multiput(8.37,35.19)(-.033714286,-.037198157){217}{\line(0,-1){.037198157}}
\multiput(1.138,24.932)(.033714286,-.044175115){217}{\line(0,-1){.044175115}}
\put(8.454,15.346){\line(0,1){0}}
\put(1.077,31.234){\circle{2.154}}
\multiput(1.09,30.376)(.033635514,-.035672897){214}{\line(0,-1){.035672897}}
\multiput(.982,32.338)(.08426136,.03346591){88}{\line(1,0){.08426136}}
\put(8.997,6.389){\makebox(0,0)[cc]{$F_{1}$}}
\put(12.049,35.112){\makebox(0,0)[cc]{$1$}}
\put(12.049,21.815){\makebox(0,0)[cc]{$a$}}
\put(12.145,15.184){\makebox(0,0)[cc]{$0$}}
\put(27.439,30.379){\circle{2.154}}
\put(27.439,22.643){\circle{2.154}}
\put(27.607,15.747){\circle{2.154}}
\put(19.198,26.175){\circle{2.154}}
\put(27.607,35.929){\circle{2.154}}
\put(27.624,34.937){\line(0,-1){3.448}}
\put(27.624,29.387){\line(0,-1){5.634}}
\put(27.624,21.567){\line(0,-1){4.625}}
\multiput(26.698,35.357)(-.033631818,-.036695455){220}{\line(0,-1){.036695455}}
\multiput(19.299,25.014)(.033632558,-.043804651){215}{\line(0,-1){.043804651}}
\put(19.59,19.551){\circle{2.154}}
\multiput(19.603,20.547)(.033641791,.047741294){201}{\line(0,1){.047741294}}
\multiput(19.603,18.584)(.075505495,-.033549451){91}{\line(1,0){.075505495}}
\put(26.569,5.883){\makebox(0,0)[cc]{$F_{2}$}}
\put(30.471,35.447){\makebox(0,0)[cc]{$1$}}
\put(30.35,15.301){\makebox(0,0)[cc]{$0$}}
\put(30.44,29.955){\makebox(0,0)[cc]{$a$}}
\put(45.136,15.004){\circle{2.154}}
\put(44.968,20.553){\circle{2.154}}
\put(44.968,25.767){\circle{2.154}}
\put(44.8,30.644){\circle{2.154}}
\put(44.968,35.521){\circle{2.154}}
\put(44.732,34.444){\line(0,-1){2.691}}
\put(44.732,24.858){\line(0,-1){3.195}}
\put(44.816,29.568){\line(0,-1){2.775}}
\put(44.816,15.693){\line(0,1){0}}
\put(44.816,19.393){\line(0,-1){3.448}}
\put(38.09,31.276){\circle{2.154}}
\put(38.259,20.848){\circle{2.154}}
\multiput(38.104,32.421)(.061489362,.033648936){94}{\line(1,0){.061489362}}
\multiput(38.104,30.349)(.0425,-.033676471){136}{\line(1,0){.0425}}
\multiput(43.884,25.769)(-.049401709,-.033547009){117}{\line(-1,0){.049401709}}
\multiput(38.104,19.88)(.039460526,-.033717105){152}{\line(1,0){.039460526}}
\put(47.81,25.115){\makebox(0,0)[cc]{$a$}}
\put(44.938,5.98){\makebox(0,0)[cc]{$F_{3}$}}
\put(47.663,35.217){\makebox(0,0)[cc]{$1$}}
\put(47.978,14.962){\makebox(0,0)[cc]{$0$}}
\put(62.772,14.981){\circle{2.154}}
\put(62.604,20.531){\circle{2.154}}
\put(62.604,25.744){\circle{2.154}}
\put(62.436,30.622){\circle{2.154}}
\put(62.604,35.499){\circle{2.154}}
\put(62.368,34.422){\line(0,-1){2.691}}
\put(62.368,29.377){\line(0,-1){2.522}}
\put(62.368,24.836){\line(0,-1){3.195}}
\put(62.537,19.286){\line(0,-1){3.363}}
\put(55.204,25.912){\circle{2.154}}
\multiput(61.78,34.675)(-.033636842,-.040715789){190}{\line(0,-1){.040715789}}
\multiput(55.305,24.92)(.033723958,-.052552083){192}{\line(0,-1){.052552083}}
\put(55.319,19.736){\circle{2.154}}
\put(55.151,31.076){\circle{2.154}}
\multiput(55.055,32.072)(.08669231,.03355128){78}{\line(1,0){.08669231}}
\multiput(55.165,30.108)(.047815385,-.033546154){130}{\line(1,0){.047815385}}
\multiput(61.489,25.856)(-.041605263,-.033723684){152}{\line(-1,0){.041605263}}
\multiput(55.165,18.876)(.054516667,-.033616667){120}{\line(1,0){.054516667}}
\put(65.961,25.201){\makebox(0,0)[cc]{$a$}}
\put(63.389,5.976){\makebox(0,0)[cc]{$F_{4}$}}
\put(65.487,35.194){\makebox(0,0)[cc]{$1$}}
\put(65.949,14.052){\makebox(0,0)[cc]{$0$}}
\put(82.426,29.192){\circle{2.154}}
\put(82.426,21.455){\circle{2.154}}
\put(82.594,14.56){\circle{2.154}}
\put(74.185,24.987){\circle{2.154}}
\put(82.426,34.574){\circle{2.154}}
\put(82.19,16.007){\line(0,1){0}}
\put(82.19,16.007){\line(0,1){0}}
\put(82.442,33.413){\line(0,-1){3.027}}
\put(82.442,30.386){\line(0,1){0}}
\put(82.358,28.115){\line(0,-1){5.634}}
\put(82.442,20.379){\line(0,-1){4.793}}
\multiput(74.201,23.911)(.033714286,-.044175115){217}{\line(0,-1){.044175115}}
\put(81.517,14.325){\line(0,1){0}}
\put(74.14,30.214){\circle{2.154}}
\multiput(74.154,29.355)(.033630841,-.035668224){214}{\line(0,-1){.035668224}}
\multiput(74.045,31.318)(.08426136,.03346591){88}{\line(1,0){.08426136}}
\multiput(74.23,26.151)(.08014607,.03350562){89}{\line(1,0){.08014607}}
\put(85.082,20.43){\makebox(0,0)[cc]{$a$}}
\put(85.082,29.263){\makebox(0,0)[cc]{$b$}}
\put(82.554,6.073){\makebox(0,0)[cc]{$F_5$}}
\put(85.282,33.835){\makebox(0,0)[cc]{$1$}}
\put(85.684,13.603){\makebox(0,0)[cc]{$0$}}
\put(101.147,14.103){\circle{1.95}}
\put(101.012,19.653){\circle{1.95}}
\put(101.012,24.866){\circle{1.95}}
\put(100.878,29.744){\circle{1.95}}
\put(101.012,34.621){\circle{1.95}}
\put(91.966,25.154){\circle{1.95}}
\put(96.677,25.045){\circle{1.95}}
\multiput(91.89,26.149)(.0336869919,.0341341463){246}{\line(0,1){.0341341463}}
\put(100.177,34.546){\line(0,1){0}}
\multiput(91.89,24.186)(.0336305221,-.0416064257){249}{\line(0,-1){.0416064257}}
\multiput(96.6,26.149)(.033693069,.034554455){101}{\line(0,1){.034554455}}
\multiput(96.687,24.077)(.033557692,-.042990385){104}{\line(0,-1){.042990385}}
\put(100.974,33.775){\line(0,-1){3.027}}
\put(100.889,28.898){\line(0,-1){3.111}}
\put(100.889,24.021){\line(0,-1){3.279}}
\put(100.974,18.807){\line(0,-1){3.7}}
\put(103.379,18.831){\makebox(0,0)[cc]{$a$}}
\put(103.487,29.845){\makebox(0,0)[cc]{$b$}}
\put(101.492,5.783){\makebox(0,0)[cc]{$F_6$}}
\put(103.913,34.316){\makebox(0,0)[cc]{$1$}}
\put(104.332,13.173){\makebox(0,0)[cc]{$0$}}
\put(120.522,13.149){\circle{1.938}}
\put(120.371,18.144){\circle{1.938}}
\put(120.371,22.836){\circle{1.938}}
\put(120.22,27.225){\circle{1.938}}
\put(120.371,31.614){\circle{1.938}}
\put(120.159,30.645){\line(0,-1){2.422}}
\put(120.159,26.105){\line(0,-1){2.27}}
\put(120.159,22.018){\line(0,-1){2.876}}
\put(120.311,17.023){\line(0,-1){3.027}}
\put(120.325,36.085){\circle{1.938}}
\put(120.29,35.261){\line(0,-1){2.801}}
\put(112.973,31.211){\circle{1.938}}
\put(112.973,17.789){\circle{1.938}}
\multiput(112.82,32.342)(.057828829,.033648649){111}{\line(1,0){.057828829}}
\multiput(112.82,30.242)(.07378161,-.03355172){87}{\line(1,0){.07378161}}
\multiput(119.356,22.656)(-.056381356,-.033635593){118}{\line(-1,0){.056381356}}
\multiput(112.82,16.937)(.057852174,-.033495652){115}{\line(1,0){.057852174}}
\put(123.235,22.429){\makebox(0,0)[cc]{$a$}}
\put(123.126,27.446){\makebox(0,0)[cc]{$b$}}
\put(120.926,5.976){\makebox(0,0)[cc]{$F_7$}}
\put(122.935,35.19){\makebox(0,0)[cc]{$1$}}
\put(123.268,13.014){\makebox(0,0)[cc]{$0$}}
\put(61.978,0){\makebox(0,0)[cc]{Figure I}}
\end{picture}

\end{center}
For $i=1,2,3,4$, let $\mathscr{B}_i(m;3,k)$ denote the subclass of $\mathscr{B}(m;3,k)$ such that $F_i$(See Figure I) is the basic block associated to $\textbf{B}\in\mathscr{B}_i(m;3,k)$. 
Following results are due to Bhavale and Aware \cite{bib2}.
\begin{proposition}\cite{bib2}\label{3redb1m}
 For an integer $m\geq 6$ and for $2\leq k\leq m-4$, $|\mathscr{B}_1(m;3,k)|=|\mathscr{B}_2(m;3,k)|=\displaystyle\sum_{l=1}^{m-5}\sum_{i=1}^{m-l-4}P^{k}_{m-l-i-2}+\sum_{r=5}^{m-2}\sum_{s=1}^{k-2}\sum_{i=1}^{r-4}P^{s+1}_{r-i-2}P^{k-s}_{m-r}.$
\end{proposition}
\begin{proposition}\cite{bib2}\label{3redb3m}
 For an integer $m\geq 7$ and for $2\leq k\leq m-5$, $\displaystyle|\mathscr{B}_3(m;3,k)|=\sum_{l=4}^{m-3}\sum_{t=1}^{k-1}P^{t+1}_{l-2}P^{k-t+1}_{m-l-1}.$
\end{proposition}
\begin{proposition}\cite{bib2}\label{3redb_4m}
For an integer $m\geq 8$ and for $3\leq k\leq m-5$, $|\mathscr{B}_{4}(m;3,k)|=\displaystyle\sum_{r=1}^{m-7}\sum_{l=4}^{m-r-3}\sum_{t=1}^{k-2}P^{t+1}_{l-2}P^{k-t}_{m-r-l-1}+\sum_{r=2}^{m-7}\sum_{s=2}^{k-2}\sum_{l=4}^{m-r-3}\sum_{t=1}^{k-s-1} P^{t+1}_{l-2}P^{k-s-t+1}_{m-r-l-1}P^{s}_{r}.$
\end{proposition}
The following result follows from Proposition \ref{3redb1m}, Proposition \ref{3redb3m}, and Proposition \ref{3redb_4m}.
\begin{theorem}\cite{bib2}\label{2redbmk}
For $m\geq 6$, $|\mathscr{B}(m;3,k)|=\displaystyle\sum_{i=1}^{4}|\mathscr{B}_i(m;3,k)|$.
\end{theorem}
The following result follows from Theorem \ref{2redbmk}.
\begin{theorem}\label{BA3red}\cite{bib2} For an integer $n\geq 6$ and for $2\leq k\leq n-4$,\\
$\displaystyle |\mathscr{L}(n;3,k)|=\displaystyle\sum_{i=0}^{n-6}(i+1)|\mathscr{B}(n-i;3,k)|$.
\end{theorem}
Thus, using Theorem \ref{2red} and Theorem \ref{BA3red}, we obtain the counting of all non-isomorphic lattices on $n$ elements, containing up to three reducible elements, and having nullity $k$. We raise the problem of counting all non-isomorphic lattices on $n$ elements, containing $r\geq 4$ reducible elements, and having nullity $k\geq 2$. Recently, Bhavale and Aware \cite{bib3} solved this problem for $r=4$ and $k=2$. Interestingly, $\mathscr{L}(n;r,k)$ consists of the lattices which are all RC-lattices for $r=4$ and $k=2$. But there exists a non RC-lattice $L\in\mathscr{L}(n;r,k)$ for $r\geq 4$ and $k\geq 3$(see \cite{bib1}). In  Section \ref{sec3}, we obtain the solution to this problem for the class of RC-lattices for which $r=4$ and $k=3$.
\section{Counting of RC-lattices containing four reducible elements and having nullity up to three}\label{sec3}
In this section, firstly we discuss counting of all non-isomorphic lattices on $n$ elements, containing four reducible elements, and having nullity two. Secondly, we count all non-isomorphic lattices on $n$ elements, containing four comparable reducible elements, and having nullity three.
\subsection{Counting of lattices containing four reducible elements and having nullity two}
According to Bhavale and Aware \cite{bib3}, if $B$ is the basic block associated to $\textbf{B}\in\mathscr{B}(j;4,2)$ then $B\in\{F_5, F_6, F_7\}$(see Figure I). Note that $F_5, F_6, F_7$ are the basic blocks associated to members of the classes $\mathscr{B}(j;4,2,3)$, $\mathscr{B}(j;4,2,4)$, $\mathscr{B}(j;4,2,5)$ respectively. The following results are due to Bhavale and Aware \cite{bib3}.
\begin{proposition}\cite{bib3}\label{Bj423}
For $j\geq 6,|\mathscr{B}(j;4,2,3)|=\binom{j-2}{4}=\displaystyle\sum_{l=1}^{j-5}\sum_{s=1}^{j-l-4}\sum_{r=1}^{j-s-l-3}(j-s-l-r-2)$.
\end{proposition}
\begin{proposition}\cite{bib3}\label{Bj424}
For $\displaystyle j\geq 7,|\mathscr{B}(j;4,2,4)|=\sum_{i=1}^{j-6}\sum_{l=2}^{j-i-4}(l-1)P^{2}_{j-i-l-2}$.
\end{proposition}
\begin{proposition}\cite{bib3}\label{Bj425}For $j\geq 8$, $\displaystyle |\mathscr{B}(j;4,2,5)|=\sum_{m=0}^{j-8}\sum_{s=4}^{j-m-4}P^{2}_{s-2}P^{2}_{j-m-s-2}$.
\end{proposition}
The following result follows from Proposition \ref{Bj423}, Proposition \ref{Bj424}, and Proposition \ref{Bj425}.
\begin{theorem}\cite{bib3}\label{4redbmk}
For $j\geq 6$,  $|\mathscr{B}(j;4,2)|=\displaystyle\sum_{i=3}^{5}|\mathscr{B}(j;4,2,i)|$.
\end{theorem}
The following result follows from Theorem \ref{4redbmk}.
\begin{theorem}\cite{bib3}\label{ln42final}
For $n\geq 6$, $\displaystyle|\mathscr{L}(n;4,2)| = \sum_{i=0}^{n-6}(i+1)|\mathscr{B}(n-i;4,2)|$.
\end{theorem}
Thus using Theorem \ref{2red}, Theorem \ref{BA3red}, and Theorem \ref{ln42final}, Bhavale and Aware \cite{bib3} obtained the counting of all non-isomorphic lattices having nullity up to two, as for $k=2$, $2\leq r\leq 4$.
\subsection{Counting of RC-lattices containing four reducible elements and having nullity three}

\begin{center}
\unitlength 1mm 
\linethickness{0.4pt}
\ifx\plotpoint\undefined\newsavebox{\plotpoint}\fi 
\begin{picture}(104.055,27.446)(0,0)
\put(5.477,10.281){\circle{.906}}
\put(5.517,14.293){\circle{.906}}
\put(5.477,18.306){\circle{.906}}
\put(5.517,22.277){\circle{.906}}
\put(5.502,21.872){\line(0,-1){3.124}}
\put(5.502,17.889){\line(0,-1){3.124}}
\put(5.502,13.848){\line(0,-1){3.124}}
\put(2.488,18.285){\circle{.906}}
\put(2.517,14.216){\circle{.906}}
\put(.453,14.244){\circle{.906}}
\multiput(5.072,22.302)(-.03349351,-.0468961){77}{\line(0,-1){.0468961}}
\multiput(2.465,17.86)(.03342308,-.04517949){78}{\line(0,-1){.04517949}}
\multiput(5.015,18.29)(-.042317757,-.033738318){107}{\line(-1,0){.042317757}}
\multiput(4.987,18.262)(-.03331081,-.0487973){74}{\line(0,-1){.0487973}}
\multiput(.43,13.791)(.042722222,-.033435185){108}{\line(1,0){.042722222}}
\multiput(2.522,13.763)(.03331081,-.04802703){74}{\line(0,-1){.04802703}}
\put(4.993,4.835){\makebox(0,0)[cc]{$B_1$}}
\put(8.094,9.985){\makebox(0,0)[cc]{$0$}}
\put(8.041,13.979){\makebox(0,0)[cc]{$a$}}
\put(8.094,18.079){\makebox(0,0)[cc]{$b$}}
\put(8.041,22.02){\makebox(0,0)[cc]{$1$}}
\put(23.994,10.782){\circle{.906}}
\put(24.034,14.794){\circle{.906}}
\put(23.994,18.806){\circle{.906}}
\put(24.034,22.778){\circle{.906}}
\put(24.019,22.373){\line(0,-1){3.124}}
\put(24.019,18.39){\line(0,-1){3.124}}
\put(24.019,14.349){\line(0,-1){3.124}}
\put(21.005,18.785){\circle{.906}}
\put(21.033,14.716){\circle{.906}}
\multiput(23.589,22.803)(-.03349351,-.0468961){77}{\line(0,-1){.0468961}}
\multiput(20.981,18.361)(.0334359,-.04519231){78}{\line(0,-1){.04519231}}
\multiput(23.503,18.762)(-.0332973,-.04878378){74}{\line(0,-1){.04878378}}
\multiput(21.039,14.263)(.0332973,-.04801351){74}{\line(0,-1){.04801351}}
\put(19.084,18.834){\circle{.906}}
\multiput(23.646,22.88)(-.041925926,-.033435185){108}{\line(-1,0){.041925926}}
\multiput(19.061,18.381)(.042722222,-.033435185){108}{\line(1,0){.042722222}}
\put(23.672,5.335){\makebox(0,0)[cc]{$B_2=B^{*}_{1}$}}
\put(32.958,16.561){\circle{.906}}
\put(39.013,10.572){\circle{.906}}
\put(39.053,14.584){\circle{.906}}
\put(39.013,18.596){\circle{.906}}
\put(39.053,22.568){\circle{.906}}
\put(39.038,22.163){\line(0,-1){3.124}}
\put(39.038,18.18){\line(0,-1){3.124}}
\put(39.038,14.139){\line(0,-1){3.124}}
\put(36.052,14.506){\circle{.906}}
\multiput(38.522,18.552)(-.0332973,-.04878378){74}{\line(0,-1){.04878378}}
\multiput(36.058,14.053)(.0332973,-.04801351){74}{\line(0,-1){.04801351}}
\multiput(32.963,16.996)(.033802395,.033634731){167}{\line(1,0){.033802395}}
\multiput(32.963,16.165)(.033662651,-.034355422){166}{\line(0,-1){.034355422}}
\put(36.023,18.425){\circle{.906}}
\multiput(38.607,22.443)(-.03349351,-.0468961){77}{\line(0,-1){.0468961}}
\multiput(35.999,18.001)(.0334359,-.04519231){78}{\line(0,-1){.04519231}}
\put(38.391,4.94){\makebox(0,0)[cc]{$B_3$}}
\put(55.764,10.631){\circle{.906}}
\put(55.804,14.643){\circle{.906}}
\put(55.764,18.655){\circle{.906}}
\put(55.804,22.627){\circle{.906}}
\put(55.789,22.222){\line(0,-1){3.124}}
\put(55.789,18.239){\line(0,-1){3.124}}
\put(55.789,14.198){\line(0,-1){3.124}}
\put(55.779,26.566){\circle{.906}}
\put(55.764,26.161){\line(0,-1){3.124}}
\put(52.806,18.55){\circle{.906}}
\multiput(55.39,22.568)(-.03349351,-.0468961){77}{\line(0,-1){.0468961}}
\multiput(52.782,18.126)(.0334359,-.04519231){78}{\line(0,-1){.04519231}}
\put(50.837,18.591){\circle{.906}}
\multiput(55.399,22.637)(-.041925926,-.033435185){108}{\line(-1,0){.041925926}}
\multiput(50.814,18.138)(.042722222,-.033435185){108}{\line(1,0){.042722222}}
\put(48.83,18.622){\circle{.906}}
\multiput(48.786,19.098)(.033694301,.038507772){193}{\line(0,1){.038507772}}
\multiput(48.786,18.206)(.03371134,-.039268041){194}{\line(0,-1){.039268041}}
\put(58.372,10.507){\makebox(0,0)[cc]{$0$}}
\put(58.315,14.405){\makebox(0,0)[cc]{$a$}}
\put(58.372,18.474){\makebox(0,0)[cc]{$x$}}
\put(58.372,22.486){\makebox(0,0)[cc]{$b$}}
\put(58.372,26.498){\makebox(0,0)[cc]{$1$}}
\put(55.215,4.963){\makebox(0,0)[cc]{$B_4$}}
\put(71.847,10.725){\circle{.906}}
\put(71.887,14.737){\circle{.906}}
\put(71.847,18.749){\circle{.906}}
\put(71.887,22.721){\circle{.906}}
\put(71.872,22.316){\line(0,-1){3.124}}
\put(71.872,18.333){\line(0,-1){3.124}}
\put(71.872,14.292){\line(0,-1){3.124}}
\put(71.862,26.66){\circle{.906}}
\put(71.847,26.255){\line(0,-1){3.124}}
\put(68.889,18.644){\circle{.906}}
\multiput(71.473,22.662)(-.03349351,-.0468961){77}{\line(0,-1){.0468961}}
\multiput(68.865,18.22)(.0334359,-.04519231){78}{\line(0,-1){.04519231}}
\put(64.913,18.716){\circle{.906}}
\multiput(64.869,19.192)(.033694301,.038507772){193}{\line(0,1){.038507772}}
\multiput(64.869,18.3)(.03371134,-.039268041){194}{\line(0,-1){.039268041}}
\put(66.846,18.678){\circle{.906}}
\multiput(71.373,26.586)(-.033559701,-.055186567){134}{\line(0,-1){.055186567}}
\multiput(66.839,18.262)(.033585185,-.055607407){135}{\line(0,-1){.055607407}}
\put(71.386,5.077){\makebox(0,0)[cc]{$B_5$}}
\put(87.209,11.058){\circle{.906}}
\put(87.249,15.07){\circle{.906}}
\put(87.209,19.082){\circle{.906}}
\put(87.249,23.054){\circle{.906}}
\put(87.234,22.649){\line(0,-1){3.124}}
\put(87.234,18.666){\line(0,-1){3.124}}
\put(87.234,14.625){\line(0,-1){3.124}}
\put(87.224,26.993){\circle{.906}}
\put(87.209,26.588){\line(0,-1){3.124}}
\put(84.215,23.028){\circle{.906}}
\multiput(86.799,27.046)(-.03349351,-.0468961){77}{\line(0,-1){.0468961}}
\multiput(84.191,22.604)(.0334359,-.04519231){78}{\line(0,-1){.04519231}}
\put(84.215,14.852){\circle{.906}}
\multiput(86.799,18.87)(-.03349351,-.0468961){77}{\line(0,-1){.0468961}}
\multiput(84.191,14.428)(.0334359,-.04519231){78}{\line(0,-1){.04519231}}
\put(81.205,21.129){\circle{.906}}
\multiput(81.16,21.642)(.034857143,.033701863){161}{\line(1,0){.034857143}}
\multiput(81.16,20.676)(.033807229,-.033584337){166}{\line(1,0){.033807229}}
\put(86.305,4.939){\makebox(0,0)[cc]{$B_6$}}
\put(103.562,10.385){\circle{.906}}
\put(103.602,14.397){\circle{.906}}
\put(103.562,18.409){\circle{.906}}
\put(103.602,22.381){\circle{.906}}
\put(103.587,21.976){\line(0,-1){3.124}}
\put(103.587,17.993){\line(0,-1){3.124}}
\put(103.587,13.952){\line(0,-1){3.124}}
\put(103.577,26.32){\circle{.906}}
\put(103.562,25.915){\line(0,-1){3.124}}
\put(100.568,22.355){\circle{.906}}
\multiput(103.152,26.373)(-.03349351,-.0468961){77}{\line(0,-1){.0468961}}
\multiput(100.544,21.931)(.0334359,-.04519231){78}{\line(0,-1){.04519231}}
\put(100.568,14.179){\circle{.906}}
\multiput(103.152,18.197)(-.03349351,-.0468961){77}{\line(0,-1){.0468961}}
\multiput(100.544,13.755)(.0334359,-.04519231){78}{\line(0,-1){.04519231}}
\put(97.602,16.181){\circle{.906}}
\multiput(97.558,16.694)(.034850932,.033701863){161}{\line(1,0){.034850932}}
\multiput(97.558,15.728)(.033801205,-.033584337){166}{\line(1,0){.033801205}}
\put(102.494,4.729){\makebox(0,0)[cc]{$B_7=B_{6}^{*}$}}
\put(54.973,.105){\makebox(0,0)[cc]{Figure II}}
\end{picture}

\unitlength 1mm 
\linethickness{0.4pt}
\ifx\plotpoint\undefined\newsavebox{\plotpoint}\fi 
\begin{picture}(105.593,30.077)(0,0)
\put(57.076,.315){\makebox(0,0)[cc]{Figure III}}
\put(7.915,6.334){\makebox(0,0)[cc]{$B_8$}}
\put(39.607,6.873){\makebox(0,0)[cc]{$B_{10}$}}
\put(23.605,6.701){\makebox(0,0)[cc]{$B_9=B_{8}^{*}$}}
\put(56.89,6.425){\makebox(0,0)[cc]{$B_{11}=B^{*}_{10}$}}
\put(73.953,6.424){\makebox(0,0)[cc]{$B_{12}$}}
\put(88.666,6.243){\makebox(0,0)[cc]{$B_{13}$}}
\put(104.563,6.334){\makebox(0,0)[cc]{$B_{14}=B^{*}_{13}$}}
\put(8.465,12.116){\circle{.884}}
\put(8.505,15.927){\circle{.884}}
\put(8.465,19.739){\circle{.884}}
\put(8.505,23.512){\circle{.884}}
\put(8.49,23.127){\line(0,-1){2.968}}
\put(8.49,19.344){\line(0,-1){2.968}}
\put(8.49,15.505){\line(0,-1){2.968}}
\put(8.48,27.254){\circle{.884}}
\put(8.465,26.87){\line(0,-1){2.968}}
\put(2.461,21.684){\circle{.884}}
\multiput(2.416,22.171)(.036679739,.03369281){153}{\line(1,0){.036679739}}
\multiput(2.416,21.253)(.035745223,-.033732484){157}{\line(1,0){.035745223}}
\put(5.444,23.445){\circle{.884}}
\multiput(8.028,27.262)(-.03349351,-.04455844){77}{\line(0,-1){.04455844}}
\multiput(5.42,23.042)(.0334359,-.0429359){78}{\line(0,-1){.0429359}}
\put(.538,21.728){\circle{.884}}
\multiput(.519,22.168)(.050472973,.033716216){148}{\line(1,0){.050472973}}
\multiput(.471,21.343)(.033626667,-.041102222){225}{\line(0,-1){.041102222}}
\put(24.737,12.138){\circle{.884}}
\put(24.777,15.949){\circle{.884}}
\put(24.737,19.761){\circle{.884}}
\put(24.777,23.534){\circle{.884}}
\put(24.762,23.149){\line(0,-1){2.968}}
\put(24.762,19.365){\line(0,-1){2.968}}
\put(24.762,15.527){\line(0,-1){2.968}}
\put(24.752,27.276){\circle{.884}}
\put(24.737,26.891){\line(0,-1){2.968}}
\put(21.743,15.742){\circle{.884}}
\multiput(24.327,19.559)(-.03349351,-.04454545){77}{\line(0,-1){.04454545}}
\multiput(21.719,15.339)(.0334359,-.04292308){78}{\line(0,-1){.04292308}}
\put(18.777,17.644){\circle{.884}}
\multiput(18.733,18.131)(.036673203,.03369281){153}{\line(1,0){.036673203}}
\multiput(18.733,17.214)(.035738854,-.033738854){157}{\line(1,0){.035738854}}
\put(16.713,17.637){\circle{.884}}
\multiput(16.662,18.023)(.033725664,.040477876){226}{\line(0,1){.040477876}}
\multiput(16.662,17.207)(.048917197,-.033643312){157}{\line(1,0){.048917197}}
\put(39.734,12.31){\circle{.884}}
\put(39.774,16.121){\circle{.884}}
\put(39.734,19.933){\circle{.884}}
\put(39.774,23.706){\circle{.884}}
\put(39.759,23.321){\line(0,-1){2.968}}
\put(39.759,19.537){\line(0,-1){2.968}}
\put(39.759,15.699){\line(0,-1){2.968}}
\put(39.749,27.448){\circle{.884}}
\put(39.734,27.063){\line(0,-1){2.968}}
\put(36.776,19.833){\circle{.884}}
\multiput(39.36,23.65)(-.03349351,-.04454545){77}{\line(0,-1){.04454545}}
\multiput(36.752,19.43)(.0334359,-.0429359){78}{\line(0,-1){.0429359}}
\put(32.8,19.901){\circle{.884}}
\multiput(32.756,20.354)(.033694301,.036580311){193}{\line(0,1){.036580311}}
\multiput(32.756,19.506)(.03371134,-.037304124){194}{\line(0,-1){.037304124}}
\put(34.733,19.865){\circle{.884}}
\multiput(34.726,19.47)(.033585185,-.05282963){135}{\line(0,-1){.05282963}}
\multiput(34.641,20.294)(.046742268,.033484536){97}{\line(1,0){.046742268}}
\put(57.498,12.016){\circle{.884}}
\put(57.538,15.828){\circle{.884}}
\put(57.498,19.639){\circle{.884}}
\put(57.538,23.413){\circle{.884}}
\put(57.523,23.028){\line(0,-1){2.968}}
\put(57.523,19.244){\line(0,-1){2.968}}
\put(57.523,15.405){\line(0,-1){2.968}}
\put(57.513,27.155){\circle{.884}}
\put(57.498,26.77){\line(0,-1){2.968}}
\put(54.54,19.539){\circle{.884}}
\multiput(57.124,23.356)(-.03349351,-.04454545){77}{\line(0,-1){.04454545}}
\multiput(54.516,19.137)(.0334359,-.0429359){78}{\line(0,-1){.0429359}}
\put(50.564,19.608){\circle{.884}}
\multiput(50.52,20.06)(.033694301,.036580311){193}{\line(0,1){.036580311}}
\multiput(50.52,19.213)(.03371134,-.037309278){194}{\line(0,-1){.037309278}}
\put(52.497,19.572){\circle{.884}}
\multiput(52.416,20.13)(.033635036,.050510949){137}{\line(0,1){.050510949}}
\multiput(52.416,19.071)(.047505155,-.033484536){97}{\line(1,0){.047505155}}
\put(74.321,11.876){\circle{.884}}
\put(74.361,15.687){\circle{.884}}
\put(74.321,19.498){\circle{.884}}
\put(74.361,23.272){\circle{.884}}
\put(74.346,22.887){\line(0,-1){2.968}}
\put(74.346,19.103){\line(0,-1){2.968}}
\put(74.346,15.264){\line(0,-1){2.968}}
\put(74.336,27.014){\circle{.884}}
\put(74.321,26.629){\line(0,-1){2.968}}
\put(68.395,21.362){\circle{.884}}
\put(68.395,17.584){\circle{.884}}
\multiput(68.296,21.772)(.035354839,.033587097){155}{\line(1,0){.035354839}}
\multiput(68.296,21.016)(.034919753,-.033691358){162}{\line(1,0){.034919753}}
\multiput(73.864,23.2)(-.036533333,-.033593333){150}{\line(-1,0){.036533333}}
\multiput(68.207,17.237)(.037217105,-.033697368){152}{\line(1,0){.037217105}}
\put(72.309,19.432){\circle{.884}}
\multiput(73.863,23.097)(-.03321277,-.06761702){47}{\line(0,-1){.06761702}}
\multiput(73.789,15.754)(-.03321277,.06910638){47}{\line(0,1){.06910638}}
\put(89.131,10.466){\circle{.884}}
\put(89.171,14.277){\circle{.884}}
\put(89.131,18.089){\circle{.884}}
\put(89.171,21.862){\circle{.884}}
\put(89.156,21.477){\line(0,-1){2.968}}
\put(89.156,17.693){\line(0,-1){2.968}}
\put(89.156,13.855){\line(0,-1){2.968}}
\put(89.146,25.604){\circle{.884}}
\put(89.131,25.219){\line(0,-1){2.968}}
\put(89.265,29.329){\circle{.884}}
\put(89.233,28.949){\line(0,-1){2.896}}
\put(84.693,25.834){\circle{.884}}
\put(84.745,14.35){\circle{.884}}
\multiput(84.713,26.253)(.04657955,.03347727){88}{\line(1,0){.04657955}}
\multiput(84.713,25.503)(.036155963,-.033431193){109}{\line(1,0){.036155963}}
\multiput(88.654,17.965)(-.039273684,-.033642105){95}{\line(-1,0){.039273684}}
\multiput(84.765,13.92)(.037394231,-.033596154){104}{\line(1,0){.037394231}}
\put(82.538,14.301){\circle{.884}}
\multiput(82.663,14.719)(.060602041,.033632653){98}{\line(1,0){.060602041}}
\multiput(88.707,10.426)(-.058009434,.033433962){106}{\line(-1,0){.058009434}}
\put(91.332,10.121){\makebox(0,0)[cc]{$0$}}
\put(91.394,13.921){\makebox(0,0)[cc]{$x$}}
\put(91.144,17.721){\makebox(0,0)[cc]{$a$}}
\put(91.269,21.402){\makebox(0,0)[cc]{$b$}}
\put(91.269,25.321){\makebox(0,0)[cc]{$y$}}
\put(91.207,28.943){\makebox(0,0)[cc]{$1$}}
\put(105.017,10.772){\circle{.884}}
\put(105.057,14.583){\circle{.884}}
\put(105.017,18.395){\circle{.884}}
\put(105.057,22.168){\circle{.884}}
\put(105.042,21.783){\line(0,-1){2.968}}
\put(105.042,17.999){\line(0,-1){2.968}}
\put(105.042,14.16){\line(0,-1){2.968}}
\put(105.032,25.91){\circle{.884}}
\put(105.017,25.525){\line(0,-1){2.968}}
\put(105.151,29.635){\circle{.884}}
\put(105.119,29.255){\line(0,-1){2.896}}
\put(100.579,26.14){\circle{.884}}
\put(100.631,14.656){\circle{.884}}
\multiput(100.599,26.559)(.04657955,.03347727){88}{\line(1,0){.04657955}}
\multiput(100.599,25.809)(.036165138,-.033431193){109}{\line(1,0){.036165138}}
\multiput(104.541,18.271)(-.039284211,-.033642105){95}{\line(-1,0){.039284211}}
\multiput(100.651,14.226)(.037403846,-.033596154){104}{\line(1,0){.037403846}}
\put(98.372,25.999){\circle{.884}}
\multiput(104.751,29.464)(-.07086517,-.03366292){89}{\line(-1,0){.07086517}}
\multiput(98.392,25.569)(.059125,-.033605769){104}{\line(1,0){.059125}}
\end{picture}

\unitlength 1mm 
\linethickness{0.4pt}
\ifx\plotpoint\undefined\newsavebox{\plotpoint}\fi 
\begin{picture}(122.264,38.042)(0,0)
\put(61.014,2.988){\makebox(0,0)[cc]{Figure IV}}
\put(3.707,17.817){\circle{.906}}
\put(3.747,21.829){\circle{.906}}
\put(3.707,25.841){\circle{.906}}
\put(3.747,29.813){\circle{.906}}
\put(3.732,29.408){\line(0,-1){3.124}}
\put(3.732,25.425){\line(0,-1){3.124}}
\put(3.732,21.384){\line(0,-1){3.124}}
\put(3.722,33.752){\circle{.906}}
\put(3.707,33.347){\line(0,-1){3.124}}
\put(.713,29.787){\circle{.906}}
\multiput(3.297,33.805)(-.03349351,-.0468961){77}{\line(0,-1){.0468961}}
\multiput(.689,29.363)(.0334359,-.04519231){78}{\line(0,-1){.04519231}}
\put(.713,21.611){\circle{.906}}
\multiput(3.297,25.629)(-.03349351,-.0468961){77}{\line(0,-1){.0468961}}
\multiput(.689,21.187)(.0334359,-.04519231){78}{\line(0,-1){.04519231}}
\put(3.668,36.913){\line(0,-1){2.694}}
\put(.453,33.756){\circle{.906}}
\multiput(.516,34.162)(.033675,.037975){80}{\line(0,1){.037975}}
\multiput(.417,33.275)(.033483516,-.077593407){91}{\line(0,-1){.077593407}}
\put(3.768,37.221){\circle{.906}}
\put(3.303,11.385){\makebox(0,0)[cc]{$B_{15}$}}
\put(17.891,17.944){\circle{.906}}
\put(17.931,21.956){\circle{.906}}
\put(17.891,25.968){\circle{.906}}
\put(17.931,29.94){\circle{.906}}
\put(17.916,29.535){\line(0,-1){3.124}}
\put(17.916,25.552){\line(0,-1){3.124}}
\put(17.916,21.511){\line(0,-1){3.124}}
\put(17.906,33.879){\circle{.906}}
\put(17.891,33.474){\line(0,-1){3.124}}
\put(17.847,37.451){\circle{.906}}
\put(17.852,37.04){\line(0,-1){2.694}}
\put(14.429,33.832){\circle{.906}}
\multiput(17.321,37.244)(-.03361905,-.03363095){84}{\line(0,-1){.03363095}}
\multiput(14.348,33.453)(.033483516,-.040021978){91}{\line(0,-1){.040021978}}
\put(14.429,25.879){\circle{.906}}
\multiput(17.395,29.737)(-.03340449,-.03675281){89}{\line(0,-1){.03675281}}
\multiput(14.348,25.426)(.033642105,-.037557895){95}{\line(0,-1){.037557895}}
\put(14.504,21.865){\circle{.906}}
\multiput(17.618,29.439)(-.03355914,-.075913978){93}{\line(0,-1){.075913978}}
\multiput(14.422,21.412)(.03889535,-.03369767){86}{\line(1,0){.03889535}}
\put(17.156,10.937){\makebox(0,0)[cc]{$B_{16}=B^{*}_{15}$}}
\put(34.964,17.733){\circle{.906}}
\put(35.004,21.745){\circle{.906}}
\put(34.964,25.757){\circle{.906}}
\put(35.004,29.729){\circle{.906}}
\put(34.989,29.324){\line(0,-1){3.124}}
\put(34.989,25.341){\line(0,-1){3.124}}
\put(34.989,21.3){\line(0,-1){3.124}}
\put(34.979,33.668){\circle{.906}}
\put(34.964,33.263){\line(0,-1){3.124}}
\put(35.098,37.589){\circle{.906}}
\put(35.066,37.189){\line(0,-1){3.048}}
\put(30.526,33.91){\circle{.906}}
\multiput(30.546,34.351)(.044554348,.033706522){92}{\line(1,0){.044554348}}
\multiput(30.546,33.562)(.034570175,-.033649123){114}{\line(1,0){.034570175}}
\multiput(34.487,25.627)(-.03731,-.03364){100}{\line(-1,0){.03731}}
\multiput(30.598,21.369)(.035354545,-.033436364){110}{\line(1,0){.035354545}}
\put(30.6,21.747){\circle{.906}}
\put(30.675,25.76){\circle{.906}}
\multiput(34.384,29.767)(-.036451923,-.033596154){104}{\line(-1,0){.036451923}}
\multiput(30.668,25.382)(.033539823,-.067752212){113}{\line(0,-1){.067752212}}
\put(34.104,10.937){\makebox(0,0)[cc]{$B_{17}$}}
\put(51.315,16.843){\circle{.906}}
\put(51.355,20.855){\circle{.906}}
\put(51.315,24.867){\circle{.906}}
\put(51.355,28.839){\circle{.906}}
\put(51.34,28.434){\line(0,-1){3.124}}
\put(51.34,24.451){\line(0,-1){3.124}}
\put(51.34,20.41){\line(0,-1){3.124}}
\put(51.33,32.778){\circle{.906}}
\put(51.315,32.373){\line(0,-1){3.124}}
\put(51.449,36.699){\circle{.906}}
\put(51.417,36.299){\line(0,-1){3.048}}
\put(46.929,20.932){\circle{.906}}
\multiput(46.897,33.461)(.044554348,.033706522){92}{\line(1,0){.044554348}}
\multiput(46.897,32.672)(.034570175,-.033649123){114}{\line(1,0){.034570175}}
\multiput(50.838,24.737)(-.03731,-.03364){100}{\line(-1,0){.03731}}
\multiput(46.949,20.479)(.035354545,-.033436364){110}{\line(1,0){.035354545}}
\put(46.431,32.898){\circle{.906}}
\put(46.728,28.661){\circle{.906}}
\multiput(50.958,36.458)(-.033626984,-.057801587){126}{\line(0,-1){.057801587}}
\multiput(46.795,28.208)(.039307692,-.033586538){104}{\line(1,0){.039307692}}
\put(50.754,10.217){\makebox(0,0)[cc]{$B_{18}=B^{*}_{17}$}}
\put(70.475,16.523){\circle{.906}}
\put(70.515,20.535){\circle{.906}}
\put(70.475,24.547){\circle{.906}}
\put(70.515,28.519){\circle{.906}}
\put(70.5,28.114){\line(0,-1){3.124}}
\put(70.5,24.131){\line(0,-1){3.124}}
\put(70.5,20.09){\line(0,-1){3.124}}
\put(70.49,32.458){\circle{.906}}
\put(70.475,32.053){\line(0,-1){3.124}}
\put(67.481,28.493){\circle{.906}}
\multiput(70.065,32.511)(-.03349351,-.0468961){77}{\line(0,-1){.0468961}}
\multiput(67.457,28.069)(.0334359,-.04519231){78}{\line(0,-1){.04519231}}
\put(67.481,20.317){\circle{.906}}
\multiput(70.065,24.335)(-.03349351,-.0468961){77}{\line(0,-1){.0468961}}
\multiput(67.457,19.893)(.0334359,-.04519231){78}{\line(0,-1){.04519231}}
\put(70.436,35.619){\line(0,-1){2.694}}
\put(70.536,35.927){\circle{.906}}
\put(64.455,28.646){\circle{.906}}
\multiput(64.522,29.234)(.033584337,.039849398){166}{\line(0,1){.039849398}}
\multiput(64.374,28.268)(.033703488,-.06525){172}{\line(0,-1){.06525}}
\put(69.923,10.086){\makebox(0,0)[cc]{$B_{19}$}}
\put(89.827,16.405){\circle{.906}}
\put(89.867,20.417){\circle{.906}}
\put(89.827,24.429){\circle{.906}}
\put(89.867,28.401){\circle{.906}}
\put(89.852,27.996){\line(0,-1){3.124}}
\put(89.852,24.013){\line(0,-1){3.124}}
\put(89.852,19.972){\line(0,-1){3.124}}
\put(89.842,32.34){\circle{.906}}
\put(89.827,31.935){\line(0,-1){3.124}}
\put(89.783,35.912){\circle{.906}}
\put(89.788,35.501){\line(0,-1){2.694}}
\put(86.365,32.293){\circle{.906}}
\multiput(89.257,35.705)(-.03361905,-.03363095){84}{\line(0,-1){.03363095}}
\multiput(86.284,31.914)(.033483516,-.040021978){91}{\line(0,-1){.040021978}}
\put(86.365,24.34){\circle{.906}}
\multiput(89.331,28.198)(-.03340449,-.03675281){89}{\line(0,-1){.03675281}}
\multiput(86.284,23.887)(.033642105,-.037557895){95}{\line(0,-1){.037557895}}
\put(83.541,24.432){\circle{.906}}
\multiput(89.554,35.5)(-.033668508,-.058723757){181}{\line(0,-1){.058723757}}
\multiput(83.534,23.905)(.033554286,-.041622857){175}{\line(0,-1){.041622857}}
\put(89.023,10.12){\makebox(0,0)[cc]{$B_{20}=B^{*}_{19}$}}
\put(106.003,16.034){\circle{.906}}
\put(106.043,20.046){\circle{.906}}
\put(106.003,24.058){\circle{.906}}
\put(106.043,28.03){\circle{.906}}
\put(106.028,27.625){\line(0,-1){3.124}}
\put(106.028,23.642){\line(0,-1){3.124}}
\put(106.028,19.601){\line(0,-1){3.124}}
\put(106.018,31.969){\circle{.906}}
\put(106.003,31.564){\line(0,-1){3.124}}
\put(106.137,35.89){\circle{.906}}
\put(106.105,35.49){\line(0,-1){3.048}}
\put(101.565,32.211){\circle{.906}}
\put(101.617,20.123){\circle{.906}}
\multiput(101.585,32.652)(.044554348,.033706522){92}{\line(1,0){.044554348}}
\multiput(101.585,31.863)(.034578947,-.033649123){114}{\line(1,0){.034578947}}
\multiput(105.527,23.928)(-.03732,-.03364){100}{\line(-1,0){.03732}}
\multiput(101.637,19.67)(.035363636,-.033436364){110}{\line(1,0){.035363636}}
\put(101.71,25.567){\circle{.906}}
\put(105.715,35.732){\line(-2,-5){3.891}}
\multiput(101.662,25.195)(.033661017,-.079677966){118}{\line(0,-1){.079677966}}
\put(119.128,15.398){\circle{.906}}
\put(119.094,18.601){\circle{.906}}
\put(119.094,21.581){\circle{.906}}
\put(119.167,24.858){\circle{.906}}
\put(119.167,28.182){\circle{.906}}
\put(119.052,31.585){\circle{.906}}
\put(119.061,35.523){\circle{.906}}
\put(122.264,34.835){\makebox(0,0)[cc]{$1$}}
\put(122.174,18.189){\makebox(0,0)[cc]{$x$}}
\put(122.078,21.274){\makebox(0,0)[cc]{$a$}}
\put(122.174,24.262){\makebox(0,0)[cc]{$y$}}
\put(122.174,28.117){\makebox(0,0)[cc]{$b$}}
\put(122.174,31.587){\makebox(0,0)[cc]{$z$}}
\put(121.982,15.214){\makebox(0,0)[cc]{$0$}}
\put(114.627,18.655){\circle{.906}}
\put(114.434,24.631){\circle{.906}}
\put(114.627,31.668){\circle{.906}}
\multiput(118.704,35.36)(-.041306122,-.03344898){98}{\line(-1,0){.041306122}}
\multiput(114.56,31.215)(.043621053,-.033484211){95}{\line(1,0){.043621053}}
\multiput(118.704,28.034)(-.04873034,-.03357303){89}{\line(-1,0){.04873034}}
\multiput(114.367,24.275)(.0542125,-.0337375){80}{\line(1,0){.0542125}}
\multiput(118.704,21.576)(-.05525333,-.03341333){75}{\line(-1,0){.05525333}}
\multiput(114.56,18.395)(.04706977,-.03362791){86}{\line(1,0){.04706977}}
\put(118.897,10.106){\makebox(0,0)[cc]{$B_{22}$}}
\put(119.088,35.118){\line(0,-1){3.135}}
\put(119.088,31.165){\line(0,-1){2.522}}
\put(119.088,27.757){\line(0,-1){2.454}}
\put(119.156,24.485){\line(0,-1){2.454}}
\put(119.088,21.146){\line(0,-1){2.181}}
\put(119.088,18.147){\line(0,-1){2.249}}
\put(105.322,9.973){\makebox(0,0)[cc]{$B_{21}$}}
\end{picture}
\end{center}
In this subsection, we count all non-isomorphic lattices on $n$ elements, containing four comparable reducible elements, and having nullity three. Let $B^{*}$ denote the dual of basic block $B$. According to Bhavale \cite{bib1}, there are exactly twenty two basic blocks (see Figure II, Figure III, and Figure IV) containing four comparable reducible elements, and having nullity three. Therefore we have the following result.
\begin{proposition}\label{bb4} If $B$ is the basic block associated to $\textbf{B}\in \mathscr{B}(j;4,3)$ where $j\geq 7$ then $B\in \{B_1,B_2,B_3,\ldots,B_{22}\}$\textnormal{(see Figure II, Figure III, and Figure IV)}.
\end{proposition}

Now out of twenty two basic blocks there are three of height three, namely $B_1,B_2,B_3$(see Figure II), nine of height four, namely $B_4$ to $B_{12}$(see Figure II and Figure III), eleven of height five, namely $B_{13}$ to $B_{21}$(see Figure III and Figure IV), and $B_{22}$(see Figure IV) is of height six. Also, for $n\geq 7$, if $L \in \mathscr{L}(n;4,3,h)$ then $3 \leq h \leq 6$. Now in this subsection, firstly we count the classes $\mathscr{B}(j;4,3,h)$ for $3 \leq h \leq 6$. Secondly, we count the classes $\mathscr{L}(n;4,3,h)$ for $3 \leq h \leq 6$, and thereby we count the class $\mathscr{L}(n;4,3)$. For $1\leq i\leq 22$, let $\mathbb{B}_{i}(j;4,3)$ be the subclass of $\mathscr{B}(j;4,3)$ containing all maximal blocks of type $\textbf{B}\in\mathscr{B}(j;4,3)$ such that $B_i$(see Figure II, Figure III, and Figure IV) is the basic block associated to $\textbf{B}$.
\begin{remark}\label{h326}
\begin{enumerate}
\item For $j\geq 7$, $\mathscr{B}(j;4,3,3)=\displaystyle\dot\cup_{i=1}^{3}\mathbb{B}_i(j;4,3)$.
\item For $j\geq 8$, $\mathscr{B}(j;4,3,4)=\displaystyle\dot\cup_{i=4}^{12}\mathbb{B}_i(j;4,3)$.
\item For $j\geq 9$, $\mathscr{B}(j;4,3,5)=\displaystyle\dot\cup_{i=13}^{21}\mathbb{B}_i(j;4,3)$.
\item For $j\geq 10$, $\mathscr{B}(j;4,3,6)=\mathbb{B}_{22}(j;4,3)$.
\item For $j\geq7$, $\mathscr{B}(j;4,3)=\displaystyle\dot\cup_{h=3}^{6}\mathscr{B}(j;4,3,h)$.
\item For $n\geq7$, $\mathscr{L}(n;4,3)=\displaystyle\dot\cup_{h=3}^{6}\mathscr{L}(n;4,3,h)$.
\end{enumerate}
\end{remark}
\subsubsection{Counting of the class $\mathscr{B}(j;4,3,3)$}
 Now, we count the the class $\mathscr{B}(j;4,3,3)$ by counting the classes $\mathbb{B}_i(j;4,3)$ for $i=1$ to $3$.
 For $x<y$, the interval $[x,y)=\{a\in P:x\leq a<y\}$, and $(x,y)=\{a\in P:x<a<y\}$; similarly, $(x,y]$ and $[x,y]$ can also be defined.
 \begin{proposition}\label{Bj4331}
For $j\geq 7$, $\displaystyle |\mathbb{B}_1(j;4,3)|=\sum_{s=1}^{j-6}\sum_{r=1}^{j-s-5} \sum_{l=2}^{j-s-r-3}(j-s-r-l-2)P_{l}^2$.
\end{proposition}
\begin{proof}
Let $\textbf{B}\in \mathbb{B}_1(j;4,3)$. Let $0<a<b<1$ be the reducible elements of $\textbf{B}$. As $B_1$(see Figure II) is the basic block associated to $\textbf{B}$, by Theorem \ref{redb}, $Red(B_1)=Red(\textbf{B})$ and $\eta(B_1)=\eta(\textbf{B})=3$. Observe that an adjunct representation of $B_1$ is given by $B_1=C]_{a}^{1}\{c_1\}]_{0}^{b}\{c_2\}]_{0}^{b}\{c_3\}$, where $C:0\prec a\prec b\prec 1$ is a $4$-chain. Also by Theorem \ref{nulladj}, $\textbf{B}$ has an adjunct representation $\textbf{B}=C_0]_{a}^{1}C_1]_{0}^{b}C_2]_{0}^{b}C_3$, where $C_0$ is a maximal chain containing all the reducible elements of $\textbf{B}$, and $C_1, C_2, C_3$ are chains. Let $s=|[0,a)\cap \textbf{B}|\geq 1$, $t=|[a,1]\cap \textbf{B}|\geq 3$, $r=|C_1|\geq 1$, $|C_2|=l_1, |C_3|=l_2$ with $1\leq l_1\leq l_2$. Let $l=l_1+l_2$.

 Now for fixed $s,r,l$, there are $t-2=j-s-r-l-2$ choices for $b$ in $\textbf{B}$ up to isomorphism. Also note that $l$ elements can be distributed over chains $C_2$ and $C_3$ in $P_l^2$ ways up to isomorphism. That is, for fixed $s,r,l$, there are $(t-2)\times P_l^2=(j-s-r-l-2)\times P_l^2$ non-isomorphic maximal blocks in $\mathbb{B}_1(j;4,3)$. Now for fixed $s,r$, $2\leq l=j-s-t-r\leq j-s-r-3$, since $t\geq 3$.  Therefore there are $\displaystyle \sum_{l=2}^{j-s-r-3}(j-s-r-l-2)\times P_{l}^2$ maximal blocks in $\mathbb{B}_1(j;4,3)$ up to isomorphism. Now for fixed $s$, $1\leq r=j-s-t-l\leq j-s-5$, since $t\geq 3$, $l\geq 2$, and there are $\displaystyle\sum_{r=1}^{j-s-5} \sum_{l=2}^{j-s-r-3}(j-s-r-l-2)P_{l}^2$ maximal blocks in $\mathbb{B}_1(j;4,3)$ up to isomorphism. Further, $1\leq s=j-t-r-l\leq j-6$, since $t\geq 3$, $r\geq 1$, $l\geq 2$, and there are $\displaystyle \sum_{s=1}^{j-6}\sum_{r=1}^{j-s-5} \sum_{l=2}^{j-s-r-3}(j-s-r-l-2)P_{l}^2$ maximal blocks in $\mathbb{B}_1(j;4,3)$ up to isomorphism.
\end{proof}
Note that $\textbf{B}\in\mathbb{B}_2(j;4,3)$ if and only if $\textbf{B*}\in\mathbb{B}_1(j;4,3)$. Therefore by Proposition \ref{Bj4331}, we have the following result.
\begin{corollary}\label{Bj4332}
For $j\geq 7$, $\displaystyle|\mathbb{B}_2(j;4,3)|= \sum_{s=1}^{j-6}\sum_{r=1}^{j-s-5} \sum_{l=2}^{j-s-r-3}(j-s-r-l-2)P_{l}^2$.
\end{corollary}

\begin{proposition}\label{Bj4333}For $j\geq 7$, $\displaystyle|\mathbb{B}_3(j;4,3)|=\displaystyle\sum_{p=1}^{j-6}\binom{j-p-2}{4}$.
\end{proposition}
\begin{proof}
Let $\textbf{B}\in \mathbb{B}_3(j;4,3)$. Let $0<a<b<1$ be the reducible elements of $\textbf{B}$. As $B_3$(see Figure II) is the basic block associated to $\textbf{B}$, by Theorem \ref{redb}, $Red(B_3)=Red(\textbf{B})$ and $\eta(B_3)=\eta(\textbf{B})=3$. Observe that an adjunct representation of $B_3$ is given by $B_3=C]_{0}^{b}\{c_1\}]_{a}^{1}\{c_2\}]_{0}^{1}\{c_3\}$, where $C:0\prec a\prec b\prec 1$ is a $4$-chain. Also by Theorem \ref{nulladj}, $\textbf{B}$ has an adjunct representation $\textbf{B}=C_0]_{0}^{b}C_1]_{a}^{1}C_2]_{0}^{1}C_3$, where $C_0$ is a maximal chain containing all the reducible elements of $\textbf{B}$, and $C_1, C_2, C_3$ are chains.

Observe that, $\textbf{B}=\textbf{B}']_{0}^{1}C_3$, where $\textbf{B}'=C_0]_{0}^{b}C_1]_{a}^{1}C_2\in\mathscr{B}(i;4,2,3)$ for $i\geq 6$ and $|C_3|=p\geq 1$ with $j=i+p\geq 7$. Suppose $\textbf{D}=\textbf{D}']_{0}^{1}C'_3$, where $\textbf{D}'=C'_0]_{0}^{b}C'_1]_{a}^{1}C'_2\in\mathscr{B}(i;4,2,3)$ for $i\geq 6$ and $|C'_3|=p\geq 1$ with $j=i+p\geq 7$. Then it is clear that $\textbf{B}\cong \textbf{D}$ if and only if $\textbf{B}'\cong \textbf{D}'$ and $C_3\cong C'_3$. To prove this suppose $\textbf{B}\cong \textbf{D}$. As $|C_3|=|C'_3|=p$, $C_3\cong C'_3$, and hence $\textbf{B}\setminus C_3\cong \textbf{D}\setminus C'_3$, {\it{that is}}, $\textbf{B}'\cong \textbf{D}'$. The converse is obvious.

 Now for fixed $p$, there are $|\mathscr{B}(j-p;4,2,3)|$ maximal blocks in $\mathbb{B}_3(j;4,3)$ up to isomorphism. Further $1\leq p=j-i\leq j-6$, since $i\geq 6$. Therefore there are $\displaystyle\sum_{p=1}^{j-6}|\mathscr{B}(j-p;4,2,3)|$ maximal blocks in $\mathbb{B}_3(j;4,3)$ up to isomorphism. By Proposition \ref{Bj423}, $|\mathscr{B}(i;4,2,3)|=\binom{i-2}{4}$. Therefore there are $\displaystyle\sum_{p=1}^{j-6}\binom{j-p-2}{4}$ maximal blocks in $\mathbb{B}_3(j;4,3)$ up to isomorphism.
\end{proof}
Using Proposition \ref{Bj4331}, Corollary \ref{Bj4332}, and Proposition \ref{Bj4333}, we have the following result.
\begin{theorem}\label{Bj433}
For $j\geq 7$, $|\mathscr{B}(j;4,3,3)|=\displaystyle\sum_{i=1}^{3}|\mathbb{B}_i(j;4,3)|=$\\$\displaystyle \sum_{s=1}^{j-6}\sum_{r=1}^{j-s-5} \sum_{l=2}^{j-s-r-3}2(j-s-r-l-2)P_{l}^2+\sum_{p=1}^{j-6}\binom{j-p-2}{4}$.
\end{theorem}
\subsubsection{Counting of the class $\mathscr{B}(j;4,3,4)$}
Here, we count the classes $\mathbb{B}_i(j;4,3)$ for $i=4$ to $12$, consequently, we count the class $\mathscr{B}(j;4,3,4)$. For that sake, let us define $\mathscr{L}'(n;2,2)$ as the subclass of $\mathscr{L}(n;2,2)$ such that any $L\in\mathscr{L}'(n;2,2)$ is of the form $L=C\oplus \textbf{B}\oplus C'$ where $\textbf{B}$ is the maximal block, and $C,C'$ are chains with $|C|\geq 1,|C'|\geq 1$. Then we have the following result.
\begin{proposition}\label{l'n22}
For $n\geq 7$, $\displaystyle |\mathscr{L}'(n;2,2)|= \sum_{i=2}^{n-5}(i-1)P^{3}_{n-i-2}$.
\end{proposition}
\begin{proof}
Let $L\in \mathscr{L}'(n;2,2)$. Then $L=C\oplus \textbf{B}\oplus C'$ where $\textbf{B}$ is the maximal block, and $C,C'$ are chains with $|C|\geq 1,|C'|\geq 1$. Let $|C|+|C'|=i\geq 2$ and $\textbf{B}\in\mathscr{B}(p;2,2)$, where $p=n-i\geq 5$, since $i\geq 2$. For fixed $i$, using Theorem \ref{2red}, by taking $k=2$, we have $|\mathscr{B}(n-i;2,2)|=P^{3}_{n-i-2}$. 
Now $i-2$ (excluding $0$ and $1$) elements can be distributed over the chains $C$ and $C'$ in $(i-2)+1=i-1$ ways up to isomorphism. Further, $2\leq i=n-p\leq n-5$, since $p\geq 5$. Therefore $\displaystyle |\mathscr{L}'(n;2,2)|=\sum_{i=2}^{n-5}(i-1)|\mathscr{B}(n-i;2,2)|=\sum_{i=2}^{n-5}(i-1)P^{3}_{n-i-2}$.
\end{proof}

\begin{proposition}\label{Bj4354}For $j\geq 8$, $|\mathbb{B}_4(j;4,3)|=\displaystyle \sum_{t=1}^{j-7}\sum_{i=2}^{j-t-5}(i-1)P^{3}_{j-t-i-2}$.
\end{proposition}
\begin{proof}
Let $\textbf{B}\in \mathbb{B}_4(j;4,3)$. Let $0<a<b<1$ be the reducible elements of $\textbf{B}$. As $B_4$(see Figure II) is the basic block associated to $\textbf{B}$, by Theorem \ref{redb}, $Red(B_4)=Red(\textbf{B})$ and $\eta(B_4)=\eta(\textbf{B})=3$. Observe that an adjunct representation of $B_4$ is given by $B_4=C]_{a}^{b}\{c_1\}]_{a}^{b}\{c_2\}]_{0}^{1}\{c_3\}$, where $C:0\prec a\prec x\prec b\prec 1$ is a $5$-chain. Also by Theorem \ref{nulladj}, $\textbf{B}$ has an adjunct representation $\textbf{B}=C_0]_{a}^{b}C_1]_{a}^{b}C_2]_{0}^{1}C_3$, where $C_0$ is a maximal chain containing all the reducible elements of $\textbf{B}$, and $C_1, C_2, C_3$ are chains.

Observe that, $\textbf{B}=L]_{0}^{1}C_3$ where $L\in \mathscr{L}'(m;2,2)$ for $m\geq 7$ and $|C_3|=t\geq 1$ with $j=m+t\geq 8$. Suppose $\textbf{D}=L']_{0}^{1}C'_3$ where $L'\in \mathscr{L}'(m;2,2)$ for $m\geq 7$ and $|C'_3|=t\geq 1$ with $j=m+t\geq 8$. Then it is clear that $\textbf{B}\cong \textbf{D}$ if and only if $L\cong L'$ and
$C_3\cong C'_3$.

 Now for fixed $t$, there are $|\mathscr{L}'(j-t;2,2)|$ maximal blocks in $\mathbb{B}_4(j;4,3)$ of type $\textbf{B}$ up to isomorphism. Further $1\leq t=j-m\leq j-7$, since $m\geq 7$. Therefore $\displaystyle |\mathbb{B}_4(j;4,3)|=\sum_{t=1}^{j-7}|\mathscr{L}'(j-t;2,2)|$. But by Proposition \ref{l'n22}, $\displaystyle |\mathscr{L}'(j-t;2,1)|= \sum_{i=2}^{j-t-5}(i-1)P^{3}_{j-t-i-2}$. Hence $\displaystyle\sum_{t=1}^{j-7}\sum_{i=2}^{j-t-5}(i-1)P^{3}_{j-t-i-2}$.
\end{proof}
\begin{proposition}\label{Bj4355}For $j\geq 8$, $|\mathbb{B}_5(j;4,3)|=\displaystyle\sum_{p=4}^{j-4}\sum_{t=1}^{j-p-3}tP_{j-p-t-1}^{2}P^2_{p-2}$.
\end{proposition}
\begin{proof}
Let $\textbf{B}\in \mathbb{B}_5(j;4,3)$. Let $0<a<b<1$ be the reducible elements of $\textbf{B}$. As $B_5$(see Figure II) is the basic block associated to $\textbf{B}$, by Theorem \ref{redb}, $Red(B_5)=Red(\textbf{B})$ and $\eta(B_5)=\eta(\textbf{B})=3$. Observe that an adjunct representation of $B_5$ is given by $B_5=C]_{a}^{b}\{c_1\}]_{0}^{1}\{c_2\}]_{0}^{1}\{c_3\}$, where $C:0\prec a\prec x\prec b\prec 1$ is a $5$-chain. Also by Theorem \ref{nulladj}, $\textbf{B}$ has an adjunct representation $\textbf{B}=C_0]_{a}^{b}C_1]_{0}^{1}C_2]_{0}^{1}C_3$, where $C_0$ is a maximal chain containing all the reducible element of $\textbf{B}$, and $C_1, C_2, C_3$ are chains. 

Observe that, $\textbf{B}=\textbf{B}']_{0}^{1}L'$ where $\textbf{B}'=C'_2]_{0}^{1}C_3\in \mathscr{B}(p;2,1)$ for $p\geq 4$ with $C'_2=\{0\}\oplus C_2\oplus\{1\}$, and $L'=C_0']_{a}^{b}C_1\in \mathscr{L}(q;2,1)$ for $q\geq 4$, with $C_0'=C_0\setminus\{0,1\}$. Note that $j=p+q\geq 8$. Suppose $\textbf{D}=\textbf{D}']_{0}^{1}L''$ where $\textbf{D}'=C'''_2]_{0}^{1}C'_3\in \mathscr{B}(p;2,1)$ for $p\geq 4$ with $C'''_2=\{0\}\oplus C''_2\oplus\{1\}$, and $L''=C_0''']_{a}^{b}C'_1\in \mathscr{L}(q;2,1)$ for $q\geq 4$, with $C_0'''=C''_0\setminus\{0,1\}$. Then it is clear that $\textbf{B}\cong \textbf{D}$ if and only if $\textbf{B}'\cong \textbf{D}'$ and $L'\cong L''$.

 Now for fixed $p$, there are $|\mathscr{B}(p;2,1)|\times|\mathscr{L}(j-p;2,1)|$ maximal blocks in $\mathbb{B}_5(j;4,3)$ up to isomorphism. Clearly for fixed $p$, $|\mathscr{B}(p;2,1)|=P^2_{p-2}$, since $\textbf{B}'\setminus \{0,1\}=C_2\dot\cup C_3$. Further $4\leq p=j-q\leq j-4$, since $q\geq 4$. Therefore there are $\displaystyle\sum_{p=4}^{j-4}|\mathscr{B}(p;2,1)|\times|\mathscr{L}(j-p;2,1)|$ maximal blocks in $\mathbb{B}_5(j;4,3)$ up to isomorphism. Also using Theorem \ref{2red}, by taking $k=1$, we have for fixed $p$, $|\mathscr{L}(j-p;2,1)|=\displaystyle\sum_{t=1}^{j-p-3}tP_{j-p-t-1}^{2}$. Hence, there are 
$\displaystyle\sum_{p=4}^{j-4}P^2_{p-2}\times\left(\sum_{t=1}^{j-p-3}tP_{j-p-t-1}^{2}\right)=\displaystyle\sum_{p=4}^{j-4}\sum_{t=1}^{j-p-3}tP^2_{p-2}P_{j-p-t-1}^{2}$ maximal blocks in $\mathbb{B}_5(j;4,3)$ up to isomorphism.
\end{proof}

\begin{proposition}\label{Bj4356}For $j\geq 8$, \\$|\mathbb{B}_6(j;4,3)|=\displaystyle\sum_{t=1}^{j-7}\sum_{r=1}^{j-t-6}\sum_{l=1}^{j-t-r-5}\sum_{i=1}^{j-t-r-l-4}P_{j-t-r-l-i-2}^{2}$.
\end{proposition}
\begin{proof}
Let $\textbf{B}\in \mathbb{B}_6(j;4,3)$. Let $0<a<b<1$ be the reducible elements of $\textbf{B}$. As $B_6$(see Figure II) is the basic block associated to $\textbf{B}$, by Theorem \ref{redb}, $Red(B_6)=Red(\textbf{B})$ and $\eta(B_6)=\eta(\textbf{B})=3$. Observe that an adjunct representation of $B_6$ is given by $B_6=C]_{a}^{1}\{c_1\}]_{b}^{1}\{c_2\}]_{0}^{b}\{c_3\}$, where $C:0\prec a\prec b\prec y\prec 1$ is a $5$-chain. Also by Theorem \ref{nulladj}, $\textbf{B}$ has an adjunct representation $\textbf{B}=C_0]_{a}^{1}C_1]_{b}^{1}C_2]_{0}^{b}C_3$, where $C_0$ is a maximal chain containing all the reducible elements of $\textbf{B}$, and $C_1,C_2,C_3$ are chains.

Observe that, $\textbf{B}=L]_{0}^{b}C_3$ where $L=C'\oplus \textbf{B}'$ with $\textbf{B}'\in\mathscr{B}_{1}(m;3,2)$ for $m\geq 6$, $C'$ is a chain with $|C'|=r\geq 1$, $C_3$ is a chain with $|C_3|=t\geq 1$, and $j=m+r+t\geq 8$. Suppose $\textbf{D}=L']_{0}^{b}C'_3$ where $L'=C''\oplus \textbf{D}'$ with $\textbf{D}'\in\mathscr{B}_{1}(m;3,2)$ for $m\geq 6$, $C''$ is a chain with $|C''|=r\geq 1$, $C'_3$ is a chain with $|C'_3|=t\geq 1$. Then it is clear that $\textbf{B}\cong \textbf{D}$ if and only if $L\cong L'$ and $C_3\cong C'_3$.

Using Proposition \ref{3redb1m}, by taking $k=2$, we have $|\mathscr{B}_{1}(m;3,2)|=\displaystyle\sum_{l=1}^{m-5}\sum_{i=1}^{m-l-4}P_{m-l-i-2}^{2}$. Therefore for fixed $r$ and $t$, there are $\displaystyle\sum_{l=1}^{j-t-r-5}\sum_{i=1}^{j-t-r-l-4}P_{j-t-r-l-i-2}^{2}$ maximal blocks of type $\textbf{B}'$ up to isomorphism in $\mathscr{B}_{1}(j-t-r;3,2)$. Now $1\leq r=j-t-m\leq j-t-6$, since $m\geq 6$. Therefore for fixed $t$, there are $\displaystyle\sum_{r=1}^{j-t-6}\sum_{l=1}^{j-t-r-5}\sum_{i=1}^{j-t-r-l-4}P_{j-t-r-l-i-2}^{2}$ lattices of type $L$ up to isomorphism. Further, $1\leq t=j-m-r\leq j-7$, since $m\geq 6$, $r\geq 1$. Therefore there are $\displaystyle\sum_{t=1}^{j-7}\sum_{r=1}^{j-t-6}\sum_{l=1}^{j-t-r-5}\sum_{i=1}^{j-t-r-l-4}P_{j-t-r-l-i-2}^{2}$ maximal blocks of type $\textbf{B}$ up to isomorphism in $\mathbb{B}_6(j;4,3)$.  Thus for $j\geq 8$, $|\mathbb{B}_6(j;4,3)|=\displaystyle\sum_{t=1}^{j-7}\sum_{r=1}^{j-t-6}\sum_{l=1}^{j-t-r-5}\sum_{i=1}^{j-t-r-l-4}P_{j-t-r-l-i-2}^{2}$.
\end{proof}
Note that $\textbf{B}\in\mathbb{B}_7(j;4,3)$ if and only if $\textbf{B*}\in\mathbb{B}_6(j;4,3)$. Therefore using Proposition \ref{Bj4356}, we have the following result.
\begin{corollary}\label{Bj4357}For $j\geq 8$, $|\mathbb{B}_7(j;4,3)|=\displaystyle\sum_{t=1}^{j-7}\sum_{r=1}^{j-t-6}\sum_{l=1}^{j-t-r-5}\sum_{i=1}^{j-t-r-l-4}P_{j-t-r-l-i-2}^{2}$.
\end{corollary}
\begin{proposition}\label{Bj4358}For $j\geq 8$, \\$|\mathbb{B}_8(j;4,3)|=\displaystyle\sum_{t=1}^{j-7}\sum_{r=1}^{j-t-6}\sum_{l=1}^{j-t-r-5}\sum_{i=1}^{j-t-r-l-4}P_{j-t-r-l-i-2}^{2}$.
\end{proposition}
\begin{proof}
Proof is similar to the proof of Proposition \ref{Bj4356}.
\end{proof}
Note that $\textbf{B}\in\mathbb{B}_9(j;4,3)$ if and only if $\textbf{B*}\in\mathbb{B}_8(j;4,3)$. Therefore using Proposition \ref{Bj4358}, we have the following result.
\begin{corollary}\label{Bj4359}For $j\geq 8$, $|\mathbb{B}_9(j;4,3)|=\displaystyle\sum_{t=1}^{j-7}\sum_{r=1}^{j-t-6}\sum_{l=1}^{j-t-r-5}\sum_{i=1}^{j-t-r-l-4}P_{j-t-r-l-i-2}^{2}$.
\end{corollary}
\begin{proposition}\label{Bj43510}For $j\geq 8$, \\$|\mathbb{B}_{10}(j;4,3)|=\displaystyle\sum_{t=1}^{j-7}\sum_{r=1}^{j-t-6}\sum_{l=1}^{j-t-r-5}\sum_{i=1}^{j-t-r-l-4}P_{j-t-r-l-i-2}^{2}$.
\end{proposition}
\begin{proof}
Proof is similar to the proof of Proposition \ref{Bj4356}.
\end{proof}
Note that $\textbf{B}\in\mathbb{B}_{11}(j;4,3)$ if and only if $\textbf{B*}\in\mathbb{B}_{10}(j;4,3)$. Therefore using Proposition \ref{Bj43510}, we have the following result.
\begin{corollary}\label{Bj43511}For $j\geq 8$, $|\mathbb{B}_{11}(j;4,3)|=\displaystyle\sum_{t=1}^{j-7}\sum_{r=1}^{j-t-6}\sum_{l=1}^{j-t-r-5}\sum_{i=1}^{j-t-r-l-4}P_{j-t-r-l-i-2}^{2}$.
\end{corollary}
\begin{proposition}\label{Bj43512}For $j\geq 8$, \\$|\mathbb{B}_{12}(j;4,3)|=\displaystyle\sum_{t=1}^{j-7}\sum_{r=1}^{j-t-6}\sum_{l=1}^{j-t-r-5}\sum_{i=1}^{j-t-r-l-4}P_{j-t-r-l-i-2}^{2}$.
\end{proposition}
\begin{proof}
Proof is similar to the proof of Proposition \ref{Bj4356}.
\end{proof}
Using Proposition \ref{Bj4354}, Proposition \ref{Bj4355}, Proposition \ref{Bj4356}, Corollary \ref{Bj4357}, Proposition \ref{Bj4358}, Corollary \ref{Bj4359}, Proposition \ref{Bj43510}, Corollary \ref{Bj43511}, and Proposition \ref{Bj43512}, we have the following result.
\begin{theorem}\label{Bj434}
For $j\geq 8$, $\mathscr{B}(j;4,3,4)=\displaystyle\sum_{i=4}^{12}|\mathbb{B}_i(j;4,3)|=\sum_{t=1}^{j-7}\sum_{i=2}^{j-t-5}(i-1)P^{3}_{j-t-i-2}$\\$+\displaystyle\sum_{p=4}^{j-4}\sum_{t=1}^{j-p-3}tP_{j-p-t-1}^{2}P^2_{p-2}+\sum_{t=1}^{j-7}\sum_{r=1}^{j-t-6}\sum_{l=1}^{j-t-r-5}\sum_{i=1}^{j-t-r-l-4}7P_{j-t-r-l-i-2}^{2}$.
\end{theorem}
\subsubsection{Counting of the class $\mathbb{B}(j;4,3,5)$}
Here, we count the classes $\mathbb{B}_i(j;4,3)$ for $i=13$ to $21$, consequently, we count the class $\mathscr{B}(j;4,3,5)$.
\begin{proposition}\label{Bj43513}
For $j\geq 9$, $\displaystyle|\mathbb{B}_{13}(j;4,3)|=\sum_{r=0}^{j-9}\sum_{p=5}^{j-r-4}P_{p-2}^{3}P_{j-p-r-2}^{2}$.
\end{proposition}
\begin{proof}
Let $\textbf{B}\in \mathbb{B}_{13}(j;4,3)$. Let $0<a<b<1$ be the reducible elements of $\textbf{B}$. As $B_{13}$(see Figure III) is the basic block associated to $\textbf{B}$, by Theorem \ref{redb}, $Red(B_{13})=Red(\textbf{B})$ and $\eta(B_{13})=\eta(\textbf{B})=3$. Observe that an adjunct representation of $B_{13}$ is given by $B_{13}=C]_{0}^{a}\{c_1\}]_{0}^{a}\{c_2\}]_{b}^{1}\{c_3\}$, where $C:0\prec x\prec a\prec b\prec y\prec 1$ is a $6$-chain. Also by Theorem \ref{nulladj}, $\textbf{B}$ has an adjunct representation $\textbf{B}=C_0]_{0}^{a}C_1]_{0}^{a}C_2]_{b}^{1}C_3$, where $C_0$ is a maximal chain containing all the reducible elements of $\textbf{B}$, and $C_1, C_2, C_3$ are chains. 

Observe that, $\textbf{B}=\textbf{B}'\oplus C'\oplus \textbf{B}''$ where $\textbf{B}'\in \mathscr{B}(p;2,2)$ with $p\geq 5$, $ \textbf{B}''\in \mathscr{B}(q;2,1)$ with $q\geq 4$, and $C'$ is a chain with $|C'|=r\geq 0$. Suppose $\textbf{D}=\textbf{D}'\oplus C''\oplus \textbf{D}''$ where $\textbf{D}'\in \mathscr{B}(p;2,2)$ with $p\geq 5$, $\textbf{D}''\in \mathscr{B}(q;2,1)$ with $q\geq 4$, and $C''$ is a chain with $|C''|=r\geq 0$. Then it is clear that $\textbf{B}\cong \textbf{D}$ if and only if $\textbf{B}'\cong \textbf{D}'$, $C'\cong C''$, and $\textbf{B}''\cong \textbf{D}''$. Note that $j=p+q+r\geq 9$. As $\textbf{B}'\in \mathscr{B}(p;2,2)$ and $\textbf{B}'\setminus\{0,a\}=(C_0\cap(0,a))\dot\cup C_1\dot\cup C_2$, $|\mathscr{B}(p;2,2)|\!=\!\!P_{p-2}^{3}$. Also as $\textbf{B}''\in \mathscr{B}(q;2,1)$ and $\textbf{B}''\setminus\{b,1\}=(C_0\cap(b,1))\dot\cup C_3$, $|\mathscr{B}(q;2,1)|=P_{q-2}^{2}$. Therefore for fixed $p$ and $r$, there are $|\mathscr{B}(p;2,2)|\times|\mathscr{B}(j-q-r;2,1)|=P_{p-2}^{3}\times P_{j-q-r-2}^{2}$ maximal blocks in $\mathbb{B}_{13}(j;4,3)$ up to isomorphism. Now $5\leq p=j-q-r\leq j-r-4$, since $q\geq 4$. Therefore for fixed $r$ by Lemma \ref{oplus}, there are $\displaystyle\sum_{p=5}^{j-r-4}|\mathscr{B}(p;2,2)|\times|\mathscr{B}(j-p-r;2,1)|=\displaystyle\sum_{p=5}^{j-r-4}P_{p-2}^{3}\times P_{j-p-r-2}^{2}$ maximal blocks in $\mathbb{B}_{13}(j;4,3)$ up to isomorphism. Further, $0\leq r=j-p-q\leq j-9$, since $p\geq 5$ and $q\geq 4$, and hence there are $\displaystyle\sum_{r=0}^{j-9}\sum_{p=5}^{j-r-4}P_{p-2}^{3}P_{j-p-r-2}^{2}$ maximal blocks in $\mathbb{B}_{13}(j;4,3)$ up to isomorphism.
\end{proof}
Note that $\textbf{B}\in\mathbb{B}_{14}(j;4,3)$ if and only if $\textbf{B*}\in\mathbb{B}_{13}(j;4,3)$. Therefore using Proposition \ref{Bj43513}, we have the following result.
\begin{corollary}\label{Bj43514}For $j\geq 9, \displaystyle|\mathbb{B}_{14}(j;4,3)|=\sum_{r=0}^{j-9}\sum_{p=5}^{j-r-4}P_{p-2}^{3}P_{j-p-r-2}^{2}$.
\end{corollary}
\begin{proposition}\label{Bj43515}
For $j\geq 9$, $\displaystyle|\mathbb{B}_{15}(j;4,3)|=\displaystyle\sum_{p=4}^{j-5}\sum_{l=1}^{j-p-4}\sum_{i=1}^{j-p-l-3}P_{p-2}^{2}P_{j-p-l-i-1}^2$.
\end{proposition}
\begin{proof}
Let $\textbf{B}\in \mathbb{B}_{15}(j;4,3)$. Let $0<a<b<1$ be the reducible elements of $\textbf{B}$. As $B_{15}$(see Figure IV) is the basic block associated to $\textbf{B}$, by Theorem \ref{redb}, $Red(B_{15})=Red(\textbf{B})$ and $\eta(B_{15})=\eta(\textbf{B})=3$. Observe that an adjunct representation of $B_{15}$ is given by $B_{15}=C]_{0}^{a}\{c_1\}]_{a}^{b}\{c_2\}]_{a}^{1}\{c_3\}$, where $C:0\prec x\prec a\prec y\prec b\prec 1$ is a $6$-chain. Also by Theorem \ref{nulladj}, $\textbf{B}$ has an adjunct representation $\textbf{B}=C_0]_{0}^{a}C_1]_{a}^{b}C_2]_{a}^{1}C_3$, where $C_0$ is a maximal chain containing all the reducible elements of $\textbf{B}$, and $C_1, C_2, C_3$ are chains.

Observe that, $\textbf{B}=\textbf{B}'\circ\textbf{B}''$ where $\textbf{B}'\in \mathscr{B}(p;2,1)$ with $p\geq 4$, $\textbf{B}''\in \mathscr{B}_{2}(q;3,2)$ with $q\geq 6$, and as the element $a$ is considered twice, $j=p+q-1$. Suppose $\textbf{D}=\textbf{D}'\circ\textbf{D}''$ where $\textbf{D}'\in \mathscr{B}(p;2,1)$ with $p\geq 4$, $\textbf{D}''\in \mathscr{B}_{2}(q;3,2)$ with $q\geq 6$. Then it is clear that $\textbf{B}\cong \textbf{D}$ if and only if $\textbf{B}'\cong \textbf{D}'$ and $\textbf{B}''\cong \textbf{D}''$. Now for fixed $p$, there are $|\mathscr{B}(p;2,1)|\times|\mathscr{B}_{2}(q;3,2)|$ maximal blocks in $\mathbb{B}_{15}(j;4,3)$ up to isomorphism, where $q=j-p+1$. Further $4\leq p=j-q+1\leq j-5$, since $q\geq 6$. Therefore by Lemma \ref{circ}, $|\mathbb{B}_{15}(j;4,3)|=\displaystyle\sum_{p=4}^{j-5}|\mathscr{B}(p;2,1)|\times|\mathscr{B}_{2}(j-p+1;3,2)|$. As $\textbf{B}'\in \mathscr{B}(p;2,1)$ and $\textbf{B}'\setminus\{0,a\}=(C_0\cap(0,a))\dot\cup C_1$, $|\mathscr{B}(p;2,1)|=P_{p-2}^{2}$. Using Proposition \ref{3redb1m}, by taking $k=2$, we have $|\mathscr{B}_{2}(j-p+1;3,2)|=\displaystyle\sum_{l=1}^{j-p-4}\sum_{i=1}^{j-p-l-3}P_{j-p-l-i-1}^2$. Therefore $|\mathbb{B}_{15}(j;4,3)|=\displaystyle\sum_{p=4}^{j-5}\left(P_{p-2}^{2}\times\sum_{l=1}^{j-p-4}\sum_{i=1}^{j-p-l-3}P_{j-p-l-i-1}^2\right)=\displaystyle\sum_{p=4}^{j-5}\sum_{l=1}^{j-p-4}\sum_{i=1}^{j-p-l-3}P_{p-2}^{2}\times P_{j-p-l-i-1}^2$. Thus there are $\displaystyle\sum_{p=4}^{j-5}\sum_{l=1}^{j-p-4}\sum_{i=1}^{j-p-l-3}P_{p-2}^{2}P_{j-p-l-i-1}^2$ maximal blocks in $\mathbb{B}_{15}(j;4,3)$ up to isomorphism.
\end{proof}
Note that $\textbf{B}\in\mathbb{B}_{16}(j;4,3)$ if and only if $\textbf{B*}\in\mathbb{B}_{15}(j;4,3)$. Therefore using Proposition \ref{Bj43515}, we have the following result.
\begin{corollary}\label{Bj43516}For $j\geq 9$, $|\mathbb{B}_{16}(j;4,3)|=\displaystyle\sum_{p=4}^{j-5}\sum_{l=1}^{j-p-4}\sum_{i=1}^{j-p-l-3}P_{p-2}^{2}P_{j-p-l-i-1}^2$.
\end{corollary}
\begin{proposition}\label{Bj43517}
For $j\geq 9$, $\displaystyle|\mathbb{B}_{17}(j;4,3)|=\displaystyle\sum_{p=4}^{j-5}\sum_{l=1}^{j-p-4}\sum_{i=1}^{j-p-l-3}P_{p-2}^{2}P_{j-p-l-i-1}^2$.
\end{proposition}
\begin{proof}
Proof is similar to the proof of Proposition \ref{Bj43515}.
\end{proof}
Note that $\textbf{B}\in\mathbb{B}_{18}(j;4,3)$ if and only if $\textbf{B*}\in\mathbb{B}_{17}(j;4,3)$. Therefore using Proposition \ref{Bj43517}, we have the following result.
\begin{corollary}\label{Bj43518}For $j\geq 9$, $|\mathbb{B}_{18}(j;4,3)|=\displaystyle\sum_{p=4}^{j-5}\sum_{l=1}^{j-p-4}\sum_{i=1}^{j-p-l-3}P_{p-2}^{2}P_{j-p-l-i-1}^2$.
\end{corollary}
\begin{proposition}\label{Bj43519}
For $j\geq 9$, $|\mathbb{B}_{19}(j;4,3)|=\displaystyle\sum_{r=1}^{j-8}\sum_{q=1}^{j-r-7}\sum_{l=4}^{j-q-r-3}P_{l-2}^2P_{j-q-r-l-1}^2$.
\end{proposition}
\begin{proof}
Let $\textbf{B}\in \mathbb{B}_{19}(j;4,3)$. Let $0<a<b<1$ be the reducible elements of $\textbf{B}$. As $B_{19}$(see Figure IV) is the basic block associated to $\textbf{B}$, by Theorem \ref{redb}, $Red(B_{19})=Red(\textbf{B})$ and $\eta(B_{19})=\eta(\textbf{B})=3$. Observe that an adjunct representation of $B_{19}$ is given by $B_{19}=C]_{0}^{a}\{c_1\}]_{a}^{b}\{c_2\}]_{0}^{1}\{c_3\}$, where $C:0\prec x\prec a\prec y \prec b\prec 1$ is a $6$-chain. Also by Theorem \ref{nulladj}, $\textbf{B}$ has an adjunct representation $\textbf{B}=C_0]_{0}^{a}C_1]_{a}^{b}C_2]_{0}^{1}C_3$, where $C_0$ is a maximal chain containing all the reducible elements of $\textbf{B}$, and $C_1, C_2, C_3$ are chains.

Observe that, $\textbf{B}=(\textbf{B}'\oplus C')]_{0}^{1}C_3$ where $\textbf{B}'\in \mathscr{B}_{3}(p;3,2)$ with $p\geq 7$, $C'$ is a chain with $|C'|=q\geq 1$, and $|C_3|=r\geq 1$. Note that $j=p+q+r\geq 9$. Suppose $\textbf{D}=(\textbf{D}'\oplus C'')]_{0}^{1}C_3'$ where $\textbf{D}'\in \mathscr{B}_{3}(p;3,2)$ with $p\geq 7$, $C''$ is a chain with $|C''|=q\geq 1$, and $|C'_3|=r\geq 1$. Then it is clear that $\textbf{B}\cong\textbf{D}$ if and only if $\textbf{B}'\cong\textbf{D}'$, $C'\cong C''$, and $C_3\cong C_3'$. Now for fixed $q$ and $r$ there are $\displaystyle|\mathscr{B}_{3}(j-q-r;3,2)|$ maximal blocks in $\mathbb{B}_{19}(j;4,3)$ up to isomorphism. Further $1\leq q=j-p-r\leq j-r-7$, since $p\geq 7$. Therefore for fixed $r$ there are $\displaystyle\sum_{q=1}^{j-r-7}|\mathscr{B}_{3}(j-q-r;3,2)|$ maximal blocks in $\mathbb{B}_{19}(j;4,3)$ up to isomorphism. Further more $1\leq r=j-p-q\leq j-8$, since $p\geq 7$, $q\geq 1$. Therefore there are $\displaystyle\sum_{r=1}^{j-8}\sum_{q=1}^{j-r-7}|\mathscr{B}_{3}(j-q-r;3,2)|$ maximal blocks in $\mathbb{B}_{19}(j;4,3)$ up to isomorphism. Using Proposition \ref{3redb3m}, by taking $k=2$, we get $|\mathscr{B}_{3}(j-q-r;3,2)|=\displaystyle\sum_{l=4}^{j-q-r-3}P_{l-2}^2P_{j-q-r-l-1}^2$. Thus there are $\displaystyle\sum_{r=1}^{j-8}\sum_{q=1}^{j-r-7}\sum_{l=4}^{j-q-r-3}P_{l-2}^2P_{j-q-r-l-1}^2$ maximal blocks in $\mathbb{B}_{19}(j;4,3)$ up to isomorphism.
\end{proof}
Note that $\textbf{B}\in\mathbb{B}_{20}(j;4,3)$ if and only if $\textbf{B*}\in\mathbb{B}_{19}(j;4,3)$. Therefore using Proposition \ref{Bj43519}, we have the following result.
\begin{corollary}\label{Bj43520}For $j\geq 9$, $\displaystyle|\mathbb{B}_{20}(j;4,3)|=\displaystyle\sum_{r=1}^{j-8}\sum_{q=1}^{j-r-7}\sum_{l=4}^{j-q-r-3}P_{l-2}^2P_{j-q-r-l-1}^2$.
\end{corollary}
\begin{proposition}\label{Bj43521}
For $j\geq 9$,\\$|\mathbb{B}_{21}(j;4,3)|=\displaystyle\sum_{t=1}^{j-8}\sum_{m=0}^{j-t-8}\sum_{s=4}^{j-t-m-4}(j-t-m-7)P^{2}_{s-2}P^{2}_{j-t-m-s-2}$.
\end{proposition}
\begin{proof}
Let $\textbf{B}\in \mathbb{B}_{21}(j;4,3)$. Let $0<a<b<1$ be the reducible elements of $\textbf{B}$. As $B_{21}$(see Figure IV) is the basic block associated to $\textbf{B}$, by Theorem \ref{redb}, $Red(B_{21})=Red(\textbf{B})$ and $\eta(B_{21})=\eta(\textbf{B})=3$. Observe that an adjunct representation of $B_{21}$ is given by $B_{21}=C]_{0}^{a}\{c_1\}]_{b}^{1}\{c_2\}]_{0}^{1}\{c_3\}$, where $C:0\prec x\prec a\prec b\prec y\prec 1$ is a $6$-chain. Also by Theorem \ref{nulladj}, $\textbf{B}$ has an adjunct representation $\textbf{B}=C_0]_{0}^{a}C_1]_{b}^{1}C_2]_{0}^{1}C_3$, where $C_0$ is a maximal chain containing all the reducible elements of $\textbf{B}$, and $C_1, C_2, C_3$ are chains.

Observe that, $\textbf{B}=\textbf{B}']_{0}^{1}C_3$ where $\textbf{B}'\in \mathscr{B}(p;4,2,5)$ with $p\geq 8$ and $C_3$ is a chain with $|C_3|=t\geq 1$. Note that $j=p+t\geq 9$. Suppose $\textbf{D}=\textbf{D}']_{0}^{1}C_3'$ where $\textbf{D}'\in \mathscr{B}(p;4,2,5)$ with $p\geq 8$ and $C_3'$ is a chain with $|C'_3|=t\geq 1$. Then it is clear that $\textbf{B}\cong \textbf{D}$ if and only if $\textbf{B}'\cong \textbf{D}'$ and $C_3\cong C_3'$. Now for fixed $t$, there are $\displaystyle|\mathscr{B}(j-t;4,2,5)|$ maximal blocks in $\mathbb{B}_{21}(j;4,3)$ up to isomorphism. By Proposition \ref{Bj425}, we have $\displaystyle |\mathscr{B}(p;4,2,5)|=\sum_{m=0}^{p-8}\sum_{s=4}^{p-m-4}(p-m-7)P^{2}_{s-2}P^{2}_{p-m-s-2}$. Further $1\leq t=j-p\leq j-8$, since $p\geq 8$. Therefore there are $\displaystyle\sum_{t=1}^{j-8}|\mathscr{B}(j-t;4,2,5)|=\sum_{t=1}^{j-8}\left(\sum_{m=0}^{j-t-8}\sum_{s=4}^{j-t-m-4}(j-t-m-7)P^{2}_{s-2}P^{2}_{j-t-m-s-2}\right)$ maximal blocks in $\mathbb{B}_{21}(j;4,3)$ up to isomorphism.
\end{proof}
Using Proposition \ref{Bj43513}, Corollary \ref{Bj43514}, Proposition \ref{Bj43515}, Corollary \ref{Bj43516}, Proposition \ref{Bj43517}, Corollary \ref{Bj43518}, Proposition \ref{Bj43519}, Corollary \ref{Bj43520}, and Proposition \ref{Bj43521}, we have the following result.
\begin{theorem}\label{Bj435}
For $j\geq 9$, \\$|\mathscr{B}(j;4,3,5)|=\displaystyle\sum_{i=13}^{21}|\mathbb{B}_i(j;4,3)|=\sum_{r=0}^{j-9}\sum_{p=5}^{j-r-4}2P_{p-2}^{3}P_{j-p-r-2}^{2}$\\$+\displaystyle\sum_{p=4}^{j-5}\sum_{l=1}^{j-p-4}\sum_{i=1}^{j-p-l-3}4P_{p-2}^{2}P_{j-p-l-i-1}^2+\displaystyle\sum_{r=1}^{j-8}\sum_{q=1}^{j-r-7}\sum_{l=4}^{j-q-r-3}2P_{l-2}^2P_{j-q-r-l-1}^2$\\$+\displaystyle\sum_{t=1}^{j-8}\sum_{m=0}^{j-t-8}\sum_{s=4}^{j-t-m-4}(j-t-m-7)P^{2}_{s-2}P^{2}_{j-t-m-s-2}$.
\end{theorem}
\subsubsection{Counting of the class $\mathscr{L}(n;4,3)$}
Now in this subsection, firstly we count the class $\mathscr{B}(j;4,3,6)$, i.e., $\mathbb{B}_{22}(j;4,3)$, secondly we count the class $\mathscr{B}(j;4,3)$, and finally we count the class $\mathscr{L}(n;4,3)$.
\begin{proposition}\label{Bj436}
For $j\geq 10$,\\$\displaystyle|\mathscr{B}(j;4,3,6)|=|\mathbb{B}_{22}(j;4,3)|=\displaystyle\sum_{p=7}^{j-3}\sum_{l=4}^{p-3}P_{j-p-1}^2P_{l-2}^2P_{p-l-1}^2$.
\end{proposition}
\begin{proof}
Let $\textbf{B}\in \mathbb{B}_{22}(j;4,3)$. Let $0<a<b<1$ be the reducible elements of $\textbf{B}$. As $B_{22}$(see Figure IV) is the basic block associated to $\textbf{B}$, by Theorem \ref{redb}, $Red(B_{22})=Red(\textbf{B})$ and $\eta(B_{22})=\eta(\textbf{B})=3$. Observe that an adjunct representation of $B_{22}$ is given by $B_{22}=C]_{0}^{a}\{c_1\}]_{a}^{b}\{c_2\}]_{b}^{1}\{c_3\}$, where $C:0\prec x\prec a\prec y\prec b\prec z\prec 1$ is a $7$-chain. Also by Theorem \ref{nulladj}, $\textbf{B}$ has an adjunct representation $\textbf{B}=C_0]_{0}^{a}C_1]_{a}^{b}C_2]_{b}^{1}C_3$, where $C_0$ is a maximal chain containing all the reducible elements of $\textbf{B}$, and $C_1, C_2, C_3$ are chains. 

Observe that, $\textbf{B}=\textbf{B}'\circ\textbf{B}''$ where $\textbf{B}'\in \mathscr{B}_{3}(p;3,2)$ with $p\geq 7$, $ \textbf{B}''\in \mathscr{B}(q;2,1)$ with $q\geq 4$, and as the element $b$ is considered twice, $j=p+q-1\geq 10$. If $\textbf{D}=\textbf{D}'\circ\textbf{D}''$ where $\textbf{D}'\in \mathscr{B}_{3}(p;3,2)$ with $p\geq 7$, $\textbf{D}''\in \mathscr{B}(q;2,1)$ with $q\geq 4$. Then $\textbf{B}\cong \textbf{D}$ if and only if $\textbf{B}'\cong \textbf{D}'$ and $\textbf{B}''\cong \textbf{D}''$. Further, $7\leq p=j-q+1\leq j-3$, since $q\geq 4$. Hence by Lemma \ref{circ}, $\mathbb{B}_{22}(j;4,3)=\displaystyle\sum_{p=7}^{j-3}|\mathscr{B}_{3}(p;3,2)|\times|\mathscr{B}(j-p+1;2,1)|$. As $\textbf{B}''\in\mathscr{B}(j-p+1;2,1)$ and $\textbf{B}''\setminus\{0,1\}=(C_0\cap(b,1))\dot\cup C_3$, $|\mathscr{B}(j-p+1;2,1)|=P_{j-p-1}^2$. Using Proposition \ref{3redb3m} by taking $k=2$, we have $|\mathscr{B}_{3}(p;3,2)|=\displaystyle\sum_{l=4}^{p-3}P_{l-2}^2\times P_{p-l-1}^2$. Therefore $\mathbb{B}_{22}(j;4,3)=\displaystyle\sum_{p=7}^{j-3}\left(P_{j-p-1}^2\times\sum_{l=4}^{p-3}P_{l-2}^2\times P_{p-l-1}^2\right)=\displaystyle\sum_{p=7}^{j-3}\sum_{l=4}^{p-3}P_{j-p-1}^2\times P_{l-2}^2\times P_{p-l-1}^2$. Thus there are $\displaystyle\sum_{p=7}^{j-3}\sum_{l=4}^{p-3}P_{p-l-1}^2P_{l-2}^2P_{p-l-1}^2$ maximal blocks in $\mathbb{B}_{22}(j;4,3)$ up to isomorphism.
\end{proof}
As $\mathscr{B}(j;4,3)=\displaystyle\dot\cup_{h=3}^{6}\mathscr{B}(j;4,3,h)$, using Theorem \ref{Bj433}, Theorem \ref{Bj434}, Theorem \ref{Bj435}, and Proposition \ref{Bj436}, we have the following result.
\begin{theorem}\label{Bj43}
For $j\geq 7$, $|\mathscr{B}(j;4,3)|=\displaystyle\sum_{h=3}^{6}|\mathscr{B}(j;4,3,h)|=\sum_{i=1}^{22}|\mathbb{B}_{i}(j;4,3)|$.
\end{theorem}

Using Theorem \ref{Bj43}, we have the following main result.
\begin{theorem}\label{Mainthm}
For $n\geq 7$, $|\mathscr{L}(n;4,3)|=\displaystyle\sum_{i=0}^{n-7}(i+1)|\mathscr{B}(n-i;4,3)|$.
\end{theorem}
\begin{proof}
Let $L\in\mathscr{L}(n;4,3)$ with $n\geq 7$. Then $L=C\oplus\textbf{B} \oplus C'$, where $C$ and $C'$ are the chains with $|C|+|C'|=i\geq 0$, and $\textbf{B}\in\mathscr{B}(j;4,3)$ with $j=n-i\geq 7$. For fixed $i$, there are $|\mathscr{B}(n-i;4,3)|$ maximal blocks up to isomorphism. Note that, there are $i+1$ ways to arrange $i$ elements on the chains $C$ and $C'$. Thus for $0\leq i=n-j\leq n-7,\displaystyle|\mathscr{L}(n;4,3)|=\sum_{i=0}^{n-7}(i+1)|\mathscr{B}(n-i;4,3)|$. Hence the result follows from Theorem \ref{Bj43}.
\end{proof}

\bibliography{sn-bibliography}

\end{document}